\documentclass[arxiv]{imsart}

\RequirePackage{amsthm,amsmath,amsfonts,amssymb,multicol,float}
\RequirePackage[authoryear]{natbib}%% uncomment this for author-year citations
\RequirePackage[colorlinks,citecolor=blue,urlcolor=blue]{hyperref}%% uncomment this for coloring bibliography citations and linked URLs
\RequirePackage{graphicx}%% uncomment this for including figures

\startlocaldefs
\theoremstyle{plain}

\newtheorem{theorem}{Theorem}[section]

\theoremstyle{remark}

\endlocaldefs

\begin{document}
\newcommand{\boldE}{\mathbf{E}}
\newcommand{\boldP}{\mathbf{P}}
\newcommand{\blue}[1]{\textcolor{blue}{#1}}

\begin{frontmatter}
%%%%%%%%%%%%%%%%%%%%%%%%%%%%%%%%%%%%%%%%%%%%%%
%%                                          %%
%% Enter the title of your article here     %%
%%                                          %%
%%%%%%%%%%%%%%%%%%%%%%%%%%%%%%%%%%%%%%%%%%%%%%
\title{On the Tightness of Graph-based Statistics}
%\title{A sample article title with some additional note\thanksref{T1}}
\runtitle{On the Tightness of Graph-based Statistics}
%\thankstext{T1}{A sample of additional note to the title.}

\begin{aug}
%%%%%%%%%%%%%%%%%%%%%%%%%%%%%%%%%%%%%%%%%%%%%%%
%% Only one address is permitted per author. %%
%% Only division, organization and e-mail is %%
%% included in the address.                  %%
%% Additional information can be included in %%
%% the Acknowledgments section if necessary. %%
%%%%%%%%%%%%%%%%%%%%%%%%%%%%%%%%%%%%%%%%%%%%%%%
\author[A]{\fnms{Lynna} \snm{Chu}\ead[label=e1]{lchu@iastate.edu}}
\and
\author[B]{\fnms{Hao} \snm{Chen}\ead[label=e2]{hxchen@ucdavis.edu}}
%%%%%%%%%%%%%%%%%%%%%%%%%%%%%%%%%%%%%%%%%%%%%%
%% Addresses                                %%
%%%%%%%%%%%%%%%%%%%%%%%%%%%%%%%%%%%%%%%%%%%%%%
\address[A]{Department of Statistics,
Iowa State University,
\printead{e1}}

\address[B]{Department of Statistics,
University of California, Davis,
\printead{e2}}

\end{aug}

\begin{abstract}
We establish tightness of graph-based stochastic processes in the space $D[0+\epsilon,1-\epsilon]$ with $\epsilon >0$ that allows for discontinuities of the first kind. The graph-based stochastic processes are based on  statistics constructed from similarity graphs. In this setting, the classic characterization of tightness is intractable, making it difficult to obtain convergence of the limiting distributions for graph-based stochastic processes. We take an alternative approach and study the behavior of the higher moments of the graph-based test statistics. We show that, under mild conditions of the graph, tightness of the stochastic process can be established by obtaining upper bounds on the graph-based statistics' higher moments. Explicit analytical expressions for these moments are provided. The results are applicable to generic graphs, including dense graphs where the number of edges can be of higher order than the number of observations.
\end{abstract}



\end{frontmatter}
%%%%%%%%%%%%%%%%%%%%%%%%%%%%%%%%%%%%%%%%%%%%%%
%% Please use \tableofcontents for articles %%
%% with 50 pages and more                   %%
%%%%%%%%%%%%%%%%%%%%%%%%%%%%%%%%%%%%%%%%%%%%%%
%\tableofcontents

%%%%%%%%%%%%%%%%%%%%%%%%%%%%%%%%%%%%%%%%%%%%%%
%%%% Main text entry area:
\section{Introduction}

Change-point detection aims to estimate and test for the presence of change-points, locations where the distribution abrupt changes, in a sequence of observations. %Applications arise in a breadth of scientific fields, ranging from neuroscience, epidemiology, social science, and climatology to name a few, where the detection of abrupt of events may be related to important scientific phenomenon, such as disease outbreak, changes in social networks, or the detection of climate events.
Research interest in change-point problems has surged in recent years and substantial contributions by the statistics community have resulted in a range of works \citep{aue2009break, zhang2010detecting, frick2014multiscale, garreau2018consistent, wang2018high, zou2020consistent, wang2021optimal}. In particular, an area of emphasis has been given to handling complex data types such as high-dimensional data or non-Euclidean data objects, including networks and images. Most change-point methods targeting complex data types are non-parametric and aim to make minimal assumptions on the underlying data generating mechanism in order to be widely applicable without restrictive assumptions (see \cite{harchaoui2009kernel, matteson2014nonparametric, shi2018consistent, dubey2020frechet} and references therein). An obstacle for non-parametric works is that theoretical guarantees can pose immense challenges. For example, fast type I error control via analytical $p$-value approximations are generally difficult to work out in the non-parametric setting. While the increasing complexity and volume of modern datasets necessitate methods that can offer fast ways to assess changes while controlling type I error, most non-parametric approaches still depend on re-sampling techniques to obtain $p$-value approximations. 

%In this paper, we focus on establishing theoretical results for an approach to change-point detection proposed in \cite{chen2015graph} and \cite{chu2019asymptotic}, which we refer collectively to as the graph-based framework. The framework is a non-parametric testing framework which utilizes test statistics constructed from similarity graphs. The similarity graph can be provided by domain knowledge or it can be generated according to some criteria, such as a minimum spanning tree (MST) or nearest neighbor graph (NN). For a brief overview of the graph-based test framework, details can be found in Section \ref{sec:review}. 
 
Recently, a graph-based framework for change-point detection was proposed in \cite{chen2015graph} and further studied in \cite{chu2019asymptotic} that aims to address the needs of modern change-point applications by offering flexibility and fast type I error control. The framework is a non-parametric approach that utilizes test statistics constructed from similarity graphs and is applicable to any data type, including multivariate and object data, as long as a similarity measure can be defined on the sample space. The similarity graph can be provided by domain knowledge or it can be generated according to some criteria, such as the minimum spanning tree or the nearest neighbor graph. This flexibility makes the approach applicable to a broad range of problems. Moreover, simulation studies and real data applications demonstrate that the approach is powerful under many settings involving high-dimensional and non-Euclidean data types \citep{chen2015graph, chu2019asymptotic}. 

The graph-based framework is also equipped with analytical $p$-value approximations for testing the significance of change-points. This extends the graph-based frameworks applicability to settings where the volume or complexity of the observations make it computationally infeasible to assess significance. A key step in obtaining these analytical $p$-value approximations is proving, under certain regularity conditions, that the stochastic processes of the graph-based test statistics converge to Gaussian processes in finite dimensional distribution (see Theorem 3.1 in \cite{chen2015graph} and Theorem 4.1 in \cite{chu2019asymptotic}). Notably, the existing theorems do not imply \textit{convergence in distribution to Gaussian processes} since tightness of the processes is not established. Tightness guarantees the existence of limit points for weak convergence and it ensures that intervals between the time points considered in the finite-dimensional distribution are well-behaved. This is essential for the type of test statistic, the maximum scan statistic, used in this framework (see (\ref{equation:scan}) below).  

In this paper, we establish tightness of the stochastic processes for graph-based test statistics under mild conditions of the graph. In terms of theoretical work, our proof provides the final piece in establishing the limiting distribution of these graph-based processes. To do so, we derive explicit expressions for higher product moments of graph-based test statistics which are obtained by studying configurations of the graph and combinatorial analysis. Importantly, our results hold for any generic graph, including dense graphs, and can be generalized to other graph-based stochastic processes to establish weak convergence. In terms of practical applications, our results provide further confidence in utilizing the asymptotic $p$-value approximations for modern data applications and the testing of change-points.% while providing a clearer understanding of the behavior of graph-based stochastic processes.  

The paper is organized as follows: Section \ref{sec:review} provides a brief overview of the graph-based framework. The main results are given in Section \ref{sec:tight} and and the proof is provided in Section \ref{sec:proof}, with additional details in the Supplementary Material. 

\section{Review of the graph-based framework} \label{sec:review}
Let $\{ \mathbf{y_i}: i = 1, \hdots, n \}$ be a data sequence indexed by time or some other meaningful ordering, where $\mathbf{y}_t$ could be a high-dimensional observation or non-Euclidean object. In the single change-point setting, there possibly exists a change-point $\tau$ such that $\mathbf{y}_t$ follows some unknown distribution for $t \le \tau$ and follows a different (unknown) distribution for $t > \tau$. Consider that each time $t$ divides the sequence of observations into two samples: those observations before time $t$ and those observations after time $t$. The graph-based framework utilizes graph-based two-sample test statistics to test whether or not these two samples are from the same distribution. By graph-based two-sample tests we refer to tests that are based on graphs with the observations $\{\textbf{y}_i\}$ as nodes. The graph, $G$, is constructed from all observations in the sequence and is usually derived from a distance or a generalized dissimilarity on the sample space, with edges in the graph connecting observations that are ``close'' in some sense. For example, $G$ could be the minimum spanning tree (MST), which is a tree connecting all observations such that the sum of the distances of edges in the tree is minimized; $G$ could also be the nearest neighbor graph (NNG) where each observation connects to its nearest neighbors. Four statistics are considered in \cite{chen2015graph} and \cite{chu2019asymptotic}. These are based on 3 quantities of the graph which we briefly discuss below. 

 For any event $x$ let $I_x$ be the indicator function that takes $1$ if $x$ is true and $0$ otherwise. We define $g_i(t)$ as an indicator function for the event that $\mathbf{y}_i$ is observed after $t$, $g_i(t) = I_{i > t}$. For an edge $e = (i,j)$, we define 
\begin{align*}
J_e(t) = \begin{cases}
0 \hspace{5mm} \text{ if } g_i(t) \neq g_j(t) \\
1 \hspace{5mm}\text{ if } g_i(t) = g_j(t) = 0, \\
2 \hspace{5mm}\text{ if } g_i(t) = g_j(t) = 1. \\
\end{cases}
\end{align*}
For any candidate value $t$ of $\tau$, the three quantities are: 
\begin{equation}
R_0(t) = \sum_{e \, \in G} I_{J_e(t) = 0}, \hspace{5mm} R_1(t) = \sum_{e \, \in G} I_{J_e(t) = 1}, \hspace{5mm} R_2(t) = \sum_{e \, \in G} I_{J_e(t) = 2}.
\end{equation}

Then $R_0(t)$ is the number of edges connecting observations before and after $t$, $R_1(t)$ is the number of edges connecting observations prior to $t$, and $R_2(t)$ is the number of edges that connect observations after $t$.

The four statistics considered are the edge-count test statistic (\ref{eq:Z0}), generalized edge-count test statistic (\ref{eq:S}), weighted edge-count test statistic (\ref{eq:Zw}), and max-type edge-count test statistic (\ref{eq:M}): 

\begin{equation}
\label{eq:Z0}
Z(t) = -\frac{R_0(t) - \boldE(R_0(t))}{\sqrt{\textbf{Var}(R_0(t))}},
\end{equation}
\begin{equation}
\label{eq:S}
S(t) = \begin{pmatrix}
R_1(t) - \boldE(R_1(t)) \\
R_2(t) - \boldE(R_2(t)) \\
\end{pmatrix}^T  \mathbf{\Sigma}^{-1}(t) \begin{pmatrix}
R_1(t) - \boldE(R_1(t)) \\
R_2(t) - \boldE(R_2(t)) \\
\end{pmatrix}. 
\end{equation}

\begin{equation} 
\label{eq:Zw} 
Z_w(t)  = \frac{R_w(t) - \boldE(R_w(t))}{\sqrt{\text{\bf Var}(R_w(t))}},
\end{equation}
with $R_w(t)  = p(t)R_1(t) + q(t)R_2(t), \quad \quad p(t) = \frac{n-t-1}{n-2}, \quad q(t) = \frac{t-1}{n-2}$,

\begin{equation} 
\label{eq:M}
M(t) = \max \left(|Z_\text{diff}(t)|,Z_w(t)\right),
\end{equation}
where $ Z_\text{diff}(t)  = \frac{R_{\text{diff}}(t)- \boldE(R_{\text{diff}}(t))}{\sqrt{\text{\bf Var}(R_{\text{diff}}(t))} }, \text{ with } R_{\text{diff}}(t) = R_1(t)-R_2(t).$

The expected value and variance of the four test statistics are computed under the permutation null distribution and their explicit expressions can be found in \cite{chen2015graph} and \cite{chu2019asymptotic}. Each of the test statistics has its own niche where it dominates; a detailed discussion can be found in \cite{chu2019asymptotic}. %As discussed in \cite{chu2019asymptotic}, $Z(t)$ does particularly well in the presence of mean change that occurs near the middle of the sequence but can suffer from low power and biased change-point estimates in the presence of more general change (for example, scale and mean change). When the change is in mean only, but happens away from the middle of the sequence, $Z_w(t)$ tends to be more powerful than $Z(t)$. For more general changes that do not necessarily target a change in mean, the generalized edge-count statistic, $S(t)$, and the max-type statistic, $M(t)$, tend to be more powerful. All four test statistics have been shown, via simulation studies, to out-compete parametric tests as the dimension of the observations increases. 

The null hypothesis of no change-point is rejected when the maximum scan statistic 
\begin{align}
\max_{n_0 \le t \le n_1} Z_0(t), \quad \max_{n_0 \le t \le n_1} Z_w(t) , \quad \max_{n_0 \le t \le n_1} S(t) ,  \quad \max_{n_0 \le t \le n_1} M(t)  \label{equation:scan} 
\end{align} 
is greater than a threshold with $n_0$ and $n_1$ being pre-specified constraints controlling where we search for the change-point. When $n$ is small, this threshold can be obtained from permutation directly. However, this becomes computationally expensive for large $n$ and instead, \cite{chen2015graph} and \cite{chu2019asymptotic} provide accurate analytical formulas to approximate the $p$-values for these scan statistics. 

\subsection{Notation} 
Let $f_n  \precsim g_n$ denote that $f_n$ is bounded above by $g_n$ (up to a constant) asymptotically and $f_n = o(g_n)$ denote that $f_n$ is dominated by $g_n$ asymptotically. We also write $f_n = O(g_n)$ to denote that $f_n$ is bounded above and below by $g_n$, asymptotically; this will also be notated as $f_n \asymp g_n$.   

\section{Tightness of basic processes} \label{sec:tight}

\subsection{Asymptotic null distributions of the basic processes}

Given the scan statistics, we reject the null hypothesis of no change-point if the scan statistic is larger than a threshold. 
Explicitly, we are interested in the following tail probabilities: $P\left(\max_{n_0 \le t \le n_1} Z(t) > b_Z \right), \\ P\left(\max_{n_0 \le t \le n_1} S(t) > b_S \right), P\left(\max_{n_0 \le t \le n_1} Z_w(t) > b_{Z_w} \right), $ and  $P\left(\max_{n_0 \le t \le n_1} M(t) > b_M \right).$

To obtain analytical approximations of these tail probabilities, \cite{chen2015graph} and \cite{chu2019asymptotic} studied the properties of the stochastic processes $\{Z(t)\}, \{S(t)\}, \{Z_w(t)\},$ and $\{M(t)\}$ under the null hypothesis. 
Based on Lemma 3.1 in \cite{chu2019asymptotic}, $S(t)$ can be expressed as $S(t) = Z^2_w(t) + Z_\text{diff}(t)$, where $Z_w(t)$ and $Z_\text{diff}(t)$ are uncorrelated. Furthermore, $Z(t)$ can be expressed as
\begin{align*} 
Z(t) = \frac{2\sigma_{R_w}}{\sqrt{4\sigma^2_{R_w}+ (p(t)-q(t))^2\sigma^2_{R_\text{diff}}}} \times  Z_w(t)  + \frac{(p(t)-q(t))\sigma_{R_\text{diff}}}{\sqrt{4\sigma^2_{R_w}+ (p(t)-q(t))^2\sigma^2_{R_\text{diff}}}} \times Z_{\text{diff}}(t),
\end{align*} 

where $\sigma^2_{R_w} = \text{\bf Var}(R_w(t))$, $\sigma^2_{R_\text{diff}} = \text{\bf Var}(R_{\text{diff}}(t))$, and $p(t)$ and $q(t)$ are defined as in (\ref{eq:Zw}).
Therefore, these stochastic processes boil down to the basic processes: $\{Z_\text{diff}(t)\},$ and  $\{Z_w(t)\}$.% and it is sufficient to obtain the limiting distributions of $\{Z_\text{diff}(t)\},$ and  $\{Z_w(t)\}$. 

In order to show that the limiting distributions of the basic processes converge to Gaussian processes, the classic approach as presented in \cite{billingsley1968convergence} is to establish: 
\begin{enumerate}
\item The convergence of $\{Z_w(\lfloor nu \rfloor): 0 < u <1 \}, \text{ and } \{Z_\text{diff}(\lfloor nu\rfloor): 0 < u <1 \}$ to multivariate Gaussian in finite dimensional distributions. \footnote{Throughout the paper, we use $\lfloor x \rfloor$ to denote the largest integer that is no larger than x.}
\item The tightness of $\{Z_w(\lfloor nu \rfloor): 0 < u <1 \} \text{ and } \{Z_\text{diff}(\lfloor nu\rfloor ): 0 < u <1 \}$.
\end{enumerate}

The first point has been proven in \cite{chen2015graph} and \cite{chu2019asymptotic}. We prove here that the second point, tightness of the graph-based stochastic processes, does indeed hold under mild conditions for the graph. 

\subsection{Main Results} 

We first state our main results and then give an outline of the proof. We use $G$ to denote both the graph and its sets of edges. Let $G_i$ be the subgraph of $G$ containing all the edges that connect to node $\mathbf{y}_i$. Then, $|G_i|$ is the number of edges in $G_i$ of the node degree of $\mathbf{y}_i$ in $G$. The these results hold for generic similarity graphs, including dense graphs. We refer to a graph as dense if the number of edges is of higher order than the number of observations, i.e. if $|G| = O(kn)$ such that $k = O(n^\alpha)$.%, $0<\alpha < 1$.

\begin{theorem} \label{thm:tight_Zw}
Under the condition that $k$ is at least $O(1)$ and $\sum_{i=1}^n|G_i|^2 = o(kn^2)$, the stochastic process  $\{Z_w( \lfloor nu \rfloor): 0 < u <1 \}$ is tight on the space $D[0+\epsilon, 1-\epsilon]$, where $\epsilon$ is a positive constant.  
\end{theorem}

\begin{theorem} \label{thm:tight_Zdiff}
Under the condition that $k$ is at least $O(1)$ and $\sum_i|G_i|^2 - \frac{4|G|^2}{n}$ is at least $O(k^2)$, the stochastic process  $\{Z_\text{diff}( \lfloor nu \rfloor ): 0 < u <1 \}$ is tight on the space $D[0+\epsilon, 1-\epsilon]$, where $\epsilon$ is a positive constant.  
\end{theorem} 

These conditions are more relaxed than the conditions in  \cite{chen2015graph} and \cite{chu2019asymptotic} when obtaining convergence in finite dimensional distributions. 

Let $D = D[0,1]$ be the space of real functions $x$ on $[0,1]$ that are right-continuous and have left-hand limits:
\begin{enumerate}
\item[(i)] For $0 \le t < 1$, $x(t+) = \lim_{s \downarrow t} x(t)$ exists and $x(t+) = x(t)$. 
\item[(ii)] For $0 \le t < 1$, $x(t-) = \lim_{s \uparrow t} x(t)$.
\end{enumerate} Functions satisfying these two properties are known as cadlag functions. A function $x$ is said to have a discontinuity of the first kind at $t$ if the left and right limits exist but differ and $x(t)$ lies between them. Any discontinuities of a cadlag function, an element of $D$, are of the first kind. Since $$\lim_{u \downarrow c} Z_w(\lfloor nu\rfloor) = Z_w(\lfloor nc \rfloor), \hspace{5mm} \lim_{u \uparrow c} Z_w(\lfloor nu\rfloor) = Z_w(\lfloor nu\rfloor),$$ 
 $$\lim_{u \downarrow c} Z_\text{diff}( \lfloor nu \rfloor) = Z_\text{diff}(\lfloor nc \rfloor) \hspace{5mm} \lim_{u \uparrow c} Z_\text{diff}( \lfloor nu \rfloor ) = Z_\text{diff}( \lfloor nu \rfloor),$$ it follows that  $Z_w(\lfloor nu \rfloor)$ and $Z_\text{diff}( \lfloor nu \rfloor)$ are right-continuous and have left-hand limits and therefore belong to the space $D$.

The classical characterization of tightness on the space $D$ is given by Theorem 13.2 in \cite{billingsley1968convergence}, a version of which is presented here: 

\textit{A sequence of stochastic processes $\{X^n(u): 0 \le u \le 1\}$ in $D$ is tight if and only if: }
\begin{enumerate}
\item[(i)] \textit{The sequence $\{X^n(u): 0 \le u \le 1\}$ is stochastically bounded in $D$,}
\item[(ii)] \textit{For each $\epsilon > 0$, $$\lim_\delta \limsup_n P(\omega^\prime(X^n, \delta) > \epsilon) = 0,$$} 
\end{enumerate} 

\textit{where $$\omega^\prime(x,\delta) = \inf_{t_i} \max_{i} \sup_{s,t \in [t_{i-1},t_i)} |x(s) - x(t)|.$$}

In general these conditions are difficult to verify, since they involve understanding the limit supreme of a sequence. We instead take an alternative approach and use the tightness criterion proposed by Kolmogorov-Chentsov (\cite{chentsov1956weak}, Theorem 1); a variant can also be found in \cite{billingsley1968convergence}. The criterion is as follows:  

\textit{A sequence of stochastic processes $X^n(u)$, $n = 1, 2, \hdots, $ right continuous with left-hand limits, is tight if there are positive constants $C, \beta, \alpha$ not depending on $n$ such that for any $0 \le u \le v \le w \le 1$, 
\[ \boldE(|X^n(v) - X^n(u)|^{2\beta}|X^n(w) - X^n(v)|^{2\beta}) \le C(w-u)^{1+\alpha}. \] }
We set $\alpha = 1, \beta = 1$ so the condition becomes:
\begin{equation}
 \boldE\left((Z^n_w(v) - Z^n_w(u))^{2} (Z^n_w(w) - Z^n_w(v))^{2}\right)  \le C_w(w-u)^2 \label{eq:condZw}
 \end{equation} 
 \begin{equation}
 \boldE\left((Z^n_\text{diff}(v) - Z^n_\text{diff}(u))^{2} (Z^n_\text{diff}(w) - Z^n_\text{diff}(v))^{2} \right) \le C_\text{diff}(w-u)^2 \label{eq:condZdiff}
 \end{equation} 
% \begin{equation}
% \boldE((Z^n(v) - Z^n(u))^{2} (Z^n(w) - Z^n(v))^{2})  \le C(w-u)^2 \label{eq:Z}
%\end{equation} 
 where the notation $Z^n_w(u) = Z_w( \lfloor nu \rfloor)$ and $Z^n_\text{diff}(u) = Z_\text{diff}( \lfloor nu \rfloor)$.  

Both inequalities automatically hold when $(w-u) \le \frac{1}{n}$ since at least one of the following is true: (i) $\lfloor nu \rfloor = \lfloor nv \rfloor$, (ii) $\lfloor nv \rfloor = \lfloor nw \rfloor$. In what follows, we focus on the case when $(w-u) > \frac{1}{n}.$

Observe that $Z^n_w(u)$ and $Z^n_\text{diff}(u)$ are not well-defined at the boundaries, when $u = 0$ or $u = 1$. We further assume that $u,v,w = O(1)$ and therefore, cannot be too close to the boundaries. As such, we establish tightness on the domain $[0 + \epsilon, 1 + \epsilon]$, where $\epsilon$ is a positive constant. The proof of this result involves obtaining explicit expressions for the $4$th moments and product moments of $Z_w$ and $Z_\text{diff}$ using combinatorial analysis. This involves determining the different graph configurations for 4 edges to be randomly selected (with replacement) from the graph and obtaining the probabilities that each configuration will occur for the graph. Focusing on the leading terms of each configuration, we show these are bounded by $C(w-u)^2$. 

\section{Proof of Theorems \ref{thm:tight_Zw} and \ref{thm:tight_Zdiff}} \label{sec:proof} 
For simplicity, let $ \lfloor nu \rfloor = r$, $ \lfloor nv \rfloor = s$, and $ \lfloor nw \rfloor= t$ and $r<s<t$. Then, expanding (\ref{eq:condZw}), we have
\begin{align*} 
\boldE ((Z^n_w(v) -  & Z^n_w(u))^{2} (Z^n_w(w) - Z^n_w(v))^{2})   = \\
&  \boldE(Z_w^2(r)Z_w^2(s)) -  2\boldE(Z_w^2(r)Z_w(s)Z_w(t)) + \boldE(Z_w^2(r)Z_w^2(t)) \\ 
& - 2\boldE(Z_w(r)Z_w^3(s)) + \boldE(Z_w^2(s)Z_w^2(t)) - 2 \boldE(Z_w(r)Z_w(s)Z_w^2(t))  \\
&  + \boldE(Z_w^4(s)) - 2\boldE(Z_w^3(s)Z_w(t)) +  4 \boldE(Z_w(r)Z_w^2(s)Z_w(t)). 
 \end{align*} 
and similarly for $ \boldE \left((Z^n_\text{diff}(v) - Z^n_\text{diff}(u))^{2} (Z^n_\text{diff}(w) - Z^n_\text{diff}(v))^{2}\right)$ (\ref{eq:condZdiff}).

For the two basic processes, the following analytical expressions are needed for $Z_w$ 
\vspace{-2mm}
\begin{multicols}{2} \noindent
\begin{align}
& \boldE(Z_w^2(r)Z_w(s)Z_w(t)), \label{eq:one} \\
& \boldE(Z_w(r)Z_w(s)Z_w^2(t)), \\
&  \boldE(Z_w(r)Z_w^2(s)Z_w(t)), \\
&  \boldE(Z_w^2(r)Z_w^2(s)), \\
& \boldE(Z_w^2(r)Z_w^2(t)), 
\end{align}
\columnbreak
\begin{align}
& \boldE(Z_w^2(s)Z_w^2(t)), \\
& \boldE(Z_w(r)Z_w^3(s)), \\
& \boldE(Z_w^3(s)Z_w(t)), \\
& \boldE(Z_w^4(s)), \label{eq:nine}
\end{align}
\end{multicols}
\vspace{-8mm}
%\vspace{-2mm}
and the following analytical expressions are needed for $Z_\text{diff}$
\begin{multicols}{2} \noindent
\begin{align}
& \boldE(Z_\text{diff}^2(r)Z_\text{diff}(s)Z_\text{diff}(t)),  \\
& \boldE(Z_\text{diff}(r)Z_\text{diff}(s)Z_\text{diff}^2(t)), \\
&  \boldE(Z_\text{diff}(r)Z_\text{diff}^2(s)Z_\text{diff}(t)), \\
&  \boldE(Z_\text{diff}^2(r)Z_\text{diff}^2(s)), \\
& \boldE(Z_\text{diff}^2(r)Z_\text{diff}^2(t)), 
\end{align}
\columnbreak
\begin{align}
& \boldE(Z_\text{diff}^2(s)Z_\text{diff}^2(t)), \\
& \boldE(Z_\text{diff}(r)Z_\text{diff}^3(s)), \\
& \boldE(Z_\text{diff}^3(s)Z_\text{diff}(t)), \\
& \boldE(Z_\text{diff}^4(s)). \label{eq:twentyfive}
\end{align}
\end{multicols}

It is straightforward to see that all the expressions can be decomposed as combinations of $R_1$ and $R_2$. 
%For example, (\ref{eq:one}) can be expressed as: 
%\begin{align*} 
%& \boldE(Z_w^2(r)Z_w(s)Z_w(t)) = \frac{\boldE\left((R_w(r) - \boldE[R_w(t)])^2(R_w(s) - \boldE[R_w(s)])(R_w(t) - \boldE[R_w(t)])\right)}{\sqrt{\text{\bf Var}(R_w(r))\text{\bf Var}(R_w(s))\text{\bf Var}(R_w(t))}}  
%\end{align*} 
%
%Since $R_w(t)$ can be expressed as a weighted sum of $R_1(t)$ and $R_2(t)$, (\ref{eq:one}) boils down to calculating the various moments of $R_1(t)$ and $R_2(t)$. For example: 
%\begin{align*} 
% \boldE(R_w(r) & R_w(s) R_w(t))  =  q^2(r)q(s)q(t)\boldE(R_1^2(r)R_1(s)R_1(t)) \\
% & + q^2(r)p(s)p(t)\boldE(R_1^2(r)R_2(s)R_2(t))\\
%& +  q^2(r)q(s)[p(t)\boldE(R_1^2(r)R_1(s)R_2(t))+q(t)\boldE(R_1^2(r)R_2(s)R_1(t))]  \\ 
%& +2q(r)p(r)q(s)[q(t)\boldE(R_1(r)R_2(r)R_1(s)R_1(t))  + p(t)\boldE(R_1(r)R_2(r)R_1(s)R_2(t))] \\
%& +2q(r)p(r)p(s)[q(t)\boldE(R_1(r)R_2(r)R_2(s)R_1(t))  +p(t)\boldE(R_1(r)R_2(r)R_2(s)R_2(t))] \\
%& + p^2(r)q(s)[q(t)\boldE(R_2^2(r)R_1(s)R_1(t)) + p(t)\boldE(R_2^2(r)R_1(s)R_2(t))]\\ 
%& +  p^2(r)p(s)[q(t)\boldE(R_2^2(r)R_2(s)R_1(t)) + p(t)\boldE(R_2^2(r)R_2(s)R_2(t))]
%\end{align*} 
%
%The remaining terms of (\ref{eq:one}) can be decomposed into a similar way. 
Since explicit expressions for the expectation, variance, and third moments of $R_w(\cdot)$, $R_\text{diff}(\cdot)$, and $R(\cdot)$ can be found in \cite{chen2015graph} and \cite{chu2019asymptotic}, the remaining unknown quantities to be derived are the product moments  of $R_1(\cdot)$ and $R_2(\cdot)$, which can be expressed as $$\boldE(R_1^a(t^\star_1)R_2^b(t^\star_2)R_1^c(t^\star_3)R_2^d(t^\star_4))$$ where $a,b,c,d = 0,1,2,3,4$ such that $a+b+c+d = 4$ and $t^\star_1, t^\star_2, t^\star_3, t^\star_4 = r,s,t$. The full list of product moments can be found in the Supplement A.   

To derive the analytical expressions for the product moments we need to: 
\begin{enumerate} 
\item Determine different configurations for 4 edges to be randomly selected (with replacement) from the graph,
\item Derive probabilities separately for each configuration. 
\end{enumerate} 

\begin{figure}[H]
  \centering
    \includegraphics[width=0.9\textwidth]{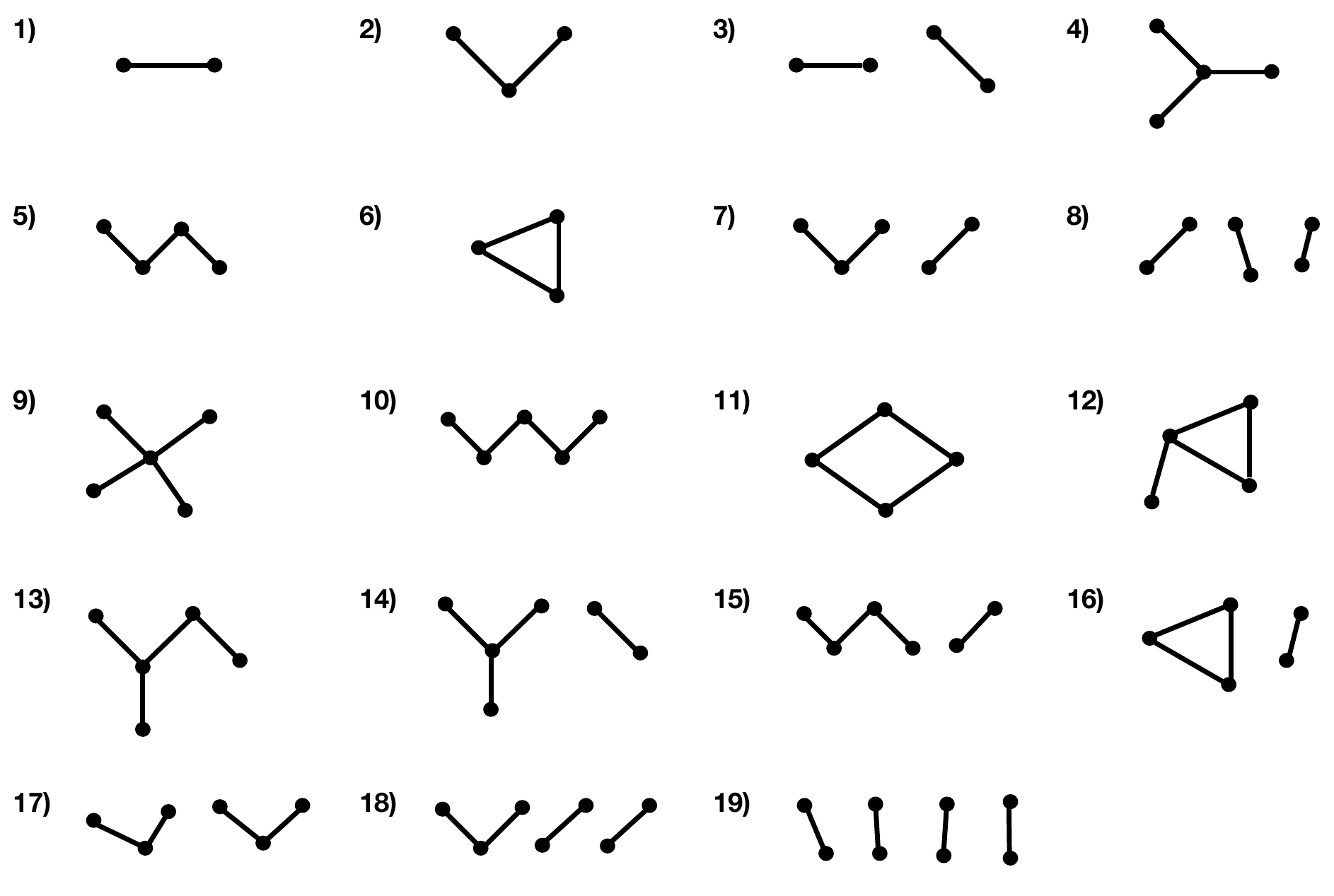}
  \caption{Nineteen configurations of 4 edges randomly chosen, with replacement, from the graph.}
  \label{fig:config}
\end{figure}

There are in total nineteen different configurations for four edges randomly chosen (with replacement) from the graph; see Figure \ref{fig:config} for an illustration of each configuration. 
Let $G$ be the similarity graph and $G_i$ be the subgraph of $G$ containing all edges that connect to node $\textbf{y}_i$. Then $|G_i|$ is the degree of node $\textbf{y}_i$ in $G$. Among all $|G|^4$ possible ways of randomly selecting the four edges, the number
of occurrences for each of the configuration are:
\begin{enumerate} 
 \item[1)]$ |G| $
 \item[2)]$ 7x_1$
 \item[3)]$  7 |G|(|G|-1) - 7x_1 $
 \item[4)]$  6x_2$
 \item[5)]$  36 x_3 $
 \item[6)]$ 12x_5$
\item[7)]$ 18x_4 - 72x_3 + 36 x_5$
\item[8)]  $ 6|G|(|G|-1)(|G|-2)- 12 x_5 - 18x_4 +  36x_3 - 6x_2$
\item[9)] $x_6$
\item[10)] $12x_7 - 24 x_8 $
\item[11)] $6x_8$
\item[12)] $24x_9$
\item[13)]$12x_{10} - 48 x_9$
\item[14)]  $4x_{11} - 12x_{10} + 24 x_9$
\item[15)] $24x_{12} - 24x_7 + 24x_8$
\item[16)] $ 8x_{13} - 24 x_9 $
\item[17)] $3x_{14} - 12x_7 + 12 x_8 $
\item[18)] $ 6x_{15} + 36x_7 - 24x_8 + 72x_9 - 12x_{10} - 48x_{12} - 24x_{13} - 6x_{14} $
\item[19)] $12x_{10} - 12x_7 - x_6 - 4x_{11} + 24x_{12} + 3x_{14} - 6x_{15} + 6x_8 + 16x_{13} +  |G|(|G|-1)(|G|-2)(|G|-3)$
 \end{enumerate}
 with $x_1, \hdots, x_{15}$ defined as: 
 \begin{align*} 
  x_1  =& \sum_{i=1}^n |G_i|^2 - 2|G|, \\
  x_2  =& \sum_{i=1}^n |G_i|^3 - 3 \sum_{i=1}^n |G_i|^2 + 4 |G|,\\
  x_3 =& \sum_{(i,j) \in G} (|G_i|-1)(|G_j|-1),\\
  x_4   =& |G| \sum_{i=1}^n |G_i|^2 + \sum_{i=1}^n |G_i|^2 - \sum_{i=1}^n |G_i|^3 - 2|G|^2, \\
  x_5 =& \sum_{(i,j)} |\{l:(i,l), (j,l) \in G \} ,  \\ 
  x_6 =& \sum_{i=1}^n |G_i|^4 - 6 \sum_{i=1}^n |G_i|^3 + 11 \sum_{i=1}^n |G_i|^2 - 12 |G|,\\
  x_7 =&  \sum_{(i,j),(j,l), i\neq l} (|G_i| -1) (|G_l|-1),  \\
  x_8 =& \sum_{(i,j), (j,l), i \neq l} |\{ m: (i,m), (l,m) \in G\}||,   \\
  x_9 =& \sum_{(i,j)} \sum_{l:(i,l), (j,l) \in G} (|G_l|-2), \\
  x_{10}  =& \sum_i \sum_{j \in G_i; j \neq i} (|G_i|-1)^2(|G_j|-1) -  2\sum_{i,j \in G} (|G_i|-1)(|G_j|-1),\\
  x_{11} =& 4|G|^2 - 3|G| \sum_{i=1}^n |G_i|^2 + |G| \sum_{i=1}^n |G_i|^3 - 2 \sum_{i=1}^n |G_i|^2 + 3 \sum_{i=1}^n |G_i|^3 - \sum_{i=1}^n |G_i|^4,\\
   x_{12}  =& |G|\sum_{(i,j)}(|G_i|-1)(|G_j|-1) - \sum_i \sum_{j \in G_i; j \neq i} (|G_i|-1)^2(|G_j|-1) \\
   & - \sum_{(i,j)}(|G_i|-1)(|G_j|-1),\\
  x_{13} =& \sum_{(i,j)}\sum_{l:(i,l), (j,l) \in G} |G \setminus \{i,j,l \in G_l\}|,  \\
x_{14}  =& \sum_i \sum_{j \neq i} (|G_i \setminus \{j \in G_i\} |) (|G_i \setminus \{j \in G_i\}|-1)(|G_j \setminus \{i \in G_j\}|)(|G_j \setminus \{i \in G_j\}|-1),  \\
  x_{15}  =& \sum_{i=1}^n  |G_i|^4 - 2 |G|\sum_{i=1}^n  |G_i|^3 + |G|^2 \sum_{i=1}^n  |G_i|^2 + |G| \sum_{i=1}^n  |G_i|^2 - \sum_{i=1}^n |G_i|^2 - 2|G|^3 + 2|G|^2. 
  \end{align*}

We will use two examples ($\boldE(R_1^2(r)R_1(s)R_1(t))$ and $\boldE(R_2(r)R_1(s)R_2(s)R_1(t)) $) to illustrate how to derive the probability for each configuration. The remaining product moments can be obtained in a similar way. The explicit formulas for all the product moments can be found in Supplement A. %\ref{supp:moments}. 
 \\

\textbf{Example 1:} To derive the probability of each configuration for $\boldE(R_1^2(r)R_1(s)R_1(t))$ (Supplement (S62)) %(\ref{eq:M1}))
, observe that 
 \begin{align*} 
& \boldE(R_1^2(r)R_1(s)R_1(t))  = \sum_{e_1,e_2,e_3,e_4} P(J_{e_1}(r) = 1, J_{e_2}(r)=1, J_{e_3}(s)=1, J_{e_4}(t)=1) \\
% & =  \hspace{-5mm} \sum_{(i,j),(l,m),(u,v),(x,y)}\hspace{-5mm}P(g_i(r) = g_j(r) = 0, g_l(r) = g_m(r) = 0, g_u(s) = g_v(s) = 0, g_x(t) = g_y(t) = 0).  
 & =  \hspace{-5mm} \sum_{\substack{(i_1,j_1),(i_2,j_2),\\(i_3,j_3),(i_4,j_4)}}\hspace{-5mm}P(g_{i_1}(r) = g_{j_1}(r) = 0, g_{i_2}(r) = g_{j_2}(r) = 0, g_{i_3}(s) = g_{j_3}(s) = 0, g_{i_4}(t) = g_{j_4}(t) = 0).  
\end{align*}
We derive $ P(g_{i_1}(r) = g_{j_1}(r) = 0, g_{i_2}(r) = g_{j_2}(r) = 0, g_{i_3}(s) = g_{j_3}(s) = 0, g_{i_4}(t) = g_{j_4}(t) = 0) \triangleq P_1$ for each of the 19 configurations separately.
\begin{enumerate}
\item[1)] The four edges are actually the same edge. 
$$P_1 = \frac{r(r-1)}{n(n-1)}.$$
\item[2)] Three edges are the same and share one node with the fourth edge or two pairs of the edges are the same and share one node.  
$$P_1 = \frac{r(r-1)((t-2)+2(s-2)+4(r-2))}{n(n-1)(n-2)}.$$
\item[3)] Three edges are the same and do not share any node with the fourth edge or two pairs of the edges are the same and do not share any node with each other.  
$$P_1 = \frac{r(r-1)((t-2)(t-3)+2(s-2)(s-3)+4(r-2)(r-3))}{n(n-1)(n-2)(n-3)}.$$
\item[4)] Two edges are the same and share one node with the other two edges. None of them share the other node (star-shaped configuration). 
$$P_1 = \frac{r(r-1)((s-2)(t-3)+2(r-2)(t-3)+3(r-2)(s-3))}{n(n-1)(n-2)(n-3)}.$$
\item[5)] Linear chain of edges such that one edge shares one node with another edge and the share the other node with the third edge. The fourth edge can be the same as any of the other three edges. 
$$P_1 = \frac{r(r-1)((r-2)(5(r-3)+4(t-3)+6(s-3))+(s-2)((s-3)+2(t-3)))}{n(n-1)(n-2)(n-3)}.$$
\item[6)] Two edges are the same and the edges form a triangle.
$$P_1 = \frac{r(r-1)((s-2)+5(r-2))}{n(n-1)(n-2)}.$$
\item[7)] Two edges share one node and do not share any node with the third edge. The fourth edge can be the same as any of the other three edges. 
\begin{align*} 
P_1 = & \frac{r(r-1)(r-2)(2(t-4)((t-3)+2(r-3))+3(s-4)((s-3)+2(r-3)))}{n(n-1)(n-2)(n-3)(n-4)}\\
& +\frac{r(r-1)(s-2)(t-4)((t-3)+2(s-3))}{n(n-1)(n-2)(n-3)(n-4)}.
\end{align*}
\item[8)] Two edges are the same and no pair of edges share any node. 
\begin{align*} 
P_1 = & \frac{r(r-1)(r-2)(r-3)(2(t-4)(t-5)+3(s-4)(s-5))}{n(n-1)(n-2)(n-3)(n-4)(n-5)}\\
&  + \frac{r(r-1)(s-2)(s-3)(t-4)(t-5)}{n(n-1)(n-2)(n-3)(n-4)(n-5)}.
\end{align*} 
\item[9)] The four edges share one node, and none of them share the other node (star-shaped). 
$$P_1 = \frac{r(r-1)(r-2)(s-3)(t-4)}{n(n-1)(n-2)(n-3)(n-4)}.$$
\item[10)] Linear chain of edges such that two distinct edges share one node with the other two edges and share a node with each other other. 
$$P_1 = \frac{r(r-1)(r-2)((r-3)((t-4)+2(s-4))+(s-3)(2(t-4)+(s-4)))}{n(n-1)(n-2)(n-3)(n-4)}.$$
\item[11)] All four edges form a box. 
$$P_1 = \frac{r(r-1)(r-2)(2(s-3)+(r-3))}{n(n-1)(n-2)(n-3)}.$$
\item[12)] Three edges form a triangle and one edge connects to one node of the triangle. 
$$P_1  = \frac{r(r-1)(r-2)(7(s-3)+2(r-3)+3(t-3))}{n(n-1)(n-2)(n-3)}.$$
\item[13)] Three edges share the same node and the fourth edge shares the other node of one of the edges. 
$$P_1 = \frac{r(r-1)(r-2)((s-3)((s-4)+7(t-4))+2(r-3)((s-4)+(t-4)))}{n(n-1)(n-2)(n-3)(n-4)}.$$
\item[14)] Three edges share the same node and the fourth edge does not share any node with the other edges. 
$$P_1 = \frac{r(r-1)(r-2)((s-3)(t-5)((s-4)+(t-4))+2(r-3)(s-4)(t-5))}{n(n-1)(n-2)(n-3)(n-4)(n-5)}.$$
\item[15)] Three edges form a linear chain and the fourth edge does not share any node with the other edges. 
\begin{align*} 
P_1  = & \frac{r(r-1)(r-2)(r-3)(3(s-4)(s-5)+4(s-4)(t-5)+(t-4)(t-5))}{n(n-1)(n-2)(n-3)(n-4)(n-5)}\\
& +\frac{2r(r-1)(r-2)(s-3)((s-4)(t-5)+(t-4)(t-5))}{n(n-1)(n-2)(n-3)(n-4)(n-5)}.
\end{align*} 
\item[16)] Three edges form a triangle and the fourth edge does not share any node with the other edges. 
$$P_1 = \frac{r(r-1)(r-2)((t-3)(t-4)+(s-3)(s-4)+2(r-3)(s-4))}{n(n-1)(n-2)(n-3)(n-4)}.$$
\item[17)] Two pairs of edges share one node with each other. The pairs of edges do not share any nodes with each other. 
$$P_1 = \frac{r(r-1)(r-2)(s-4)(t-5)((s-3)+2(r-3))}{n(n-1)(n-2)(n-3)(n-4)(n-5)}.$$
\item[18)] Two edges share one node with each other. The other edges do not share any nodes with any of the other edges.
\begin{align*} 
P_1 = & \frac{r(r-1)(r-2)(r-3)(s-4)(t-6)(2(t-5)+3(s-5))}{n(n-1)(n-2)(n-3)(n-4)(n-5)(n-6)}\\
& +\frac{r(r-1)(r-2)(s-3)(s-4)(t-5)(t-6)}{n(n-1)(n-2)(n-3)(n-4)(n-5)(n-6)}.
\end{align*}
\item[19)] None of the four edges share any node. 
$$P_1= \frac{r(r-1)(r-2)(r-3)(s-4)(s-5)(t-6)(t-7)}{n(n-1)(n-2)(n-3)(n-4)(n-5)(n-6)(n-7)}.$$
\end{enumerate} 

\textbf{Example 2:} To derive the probability of each configuration for $\boldE(R_2(r)R_1(s)R_2(s)R_1(t)) $ (Supplement (S84)) %(\ref{eq:N7}))
, observe that 
 \begin{align*} 
& \boldE(R_2(r)R_1(s)R_2(s)R_1(t))  = \sum_{e_1,e_2,e_3,e_4} P(J_{e_1}(r) = 2, J_{e_2}(s)=1, J_{e_3}(s)=2, J_{e_4}(t)=1) \\
 %& = \hspace{-5mm} \sum_{(i,j),(l,m),(u,v),(x,y)} \hspace{-7mm}P(g_i(r) = g_j(r) = 1, g_l(r) = g_m(r) = 0, g_u(s) = g_v(s) = 1, g_x(t) = g_y(t) = 0).  \\
 & =  \hspace{-5mm} \sum_{\substack{(i_1,j_1),(i_2,j_2),\\(i_3,j_3),(i_4,j_4)}}\hspace{-5mm}P(g_{i_1}(r) = g_{j_1}(r) = 1, g_{i_2}(r) = g_{j_2}(r) = 0, g_{i_3}(s) = g_{j_3}(s) = 1, g_{i_4}(t) = g_{j_4}(t) = 0).  
\end{align*}
We derive $ P(g_{i_1}(r) = g_{j_1}(r) = 1, g_{i_2}(r) = g_{j_2}(r) = 0, g_{i_3}(s) = g_{j_3}(s) = 1, g_{i_4}(t) = g_{j_4}(t) = 0) \triangleq P_2$ for each of the 19 configurations separately.
\begin{enumerate}
\item[1)] The four edges are actually the same edge. 
$$P_2 = 0.$$
\item[2)] Three edges are the same and share one node with the fourth edge or two pairs of the edges are the same and share one node.  
$$P_2 = 0.$$
\item[3)] Three edges are the same and do not share any node with the fourth edge or two pairs of the edges are the same and do not share any node with each other.  
\begin{align*} 
P_2 &= \frac{(s(s-1)+(s-r)(s-r-1))((n-s)(n-s-1)+(t-s)(t-s-1))}{n(n-1)(n-2)(n-3)} 
\end{align*} 
\item[4)] Two edges are the same and share one node with the other two edges. None of them share the other node (star-shaped configuration). 
$$P_2 = 0.$$
\item[5)] Linear chain of edges such that one edge shares one node with another edge and the share the other node with the third edge. The fourth edge can be the same as any of the other three edges. 
\begin{align*} P_2 = & \frac{(t-s)(n-s-1)((s-r)(s-r-1)+s(s-1) + (s-r)(s-1))}{n(n-1)(n-2)(n-3)}  \\
& +\frac{(s-r)(n-s-1)((t-s)(s-1)+(s-1)(n-s))}{n(n-1)(n-2)(n-3)}. 
\end{align*} 
\item[6)] Two edges are the same and the edges form a triangle.
$$P_2 = 0.$$
\item[7)] Two edges share one node and do not share any node with the third edge. The fourth edge can be the same as any of the other three edges. 
\begin{align*} 
f_{7a} = & (s-r)(s-r-1)(t-s)((s-2)(n-s-1)+(t-s-1)(n-s-2)\\
& +(n-s-1)(n-s-2)) + (s-r)(n-s-1)(2(s-r-1)(s-2)(n-s)\\
& +(t-s)(s-1)(n-s-2)),\\
f_{7b} = & (t-s)((s-r)(s-1)(s-2)(n-s-1)+(t-s-1)s(s-1)(n-s-2))\\
& +s(s-1)((s-2)(n-s)(n-s-1)+(t-s)(n-s-1)(n-s-2))\\
& +(t-s)s(s-1)((s-2)(n-s-1)+(t-s-1)(n-s-2)), \\
f_{7c}  = & r(n-s)(n-s-1)((r-1)(n-r-2)+2(s-r)(n-r-3))\\
& +(s-r)(n-s)(n-s-1)(2(s-r-1)(n-r-4)+ r(n-r-3))\\
& +(t-s)(t-s-1)(s-r)(3r(n-r-3)+2(s-r-1)(n-r-4))\\
& +(t-s)(t-s-1)r(r-1)(n-r-2), \\
P_2 = & \frac{f_{7a} + f_{7b} + f_{7c}}{n(n-1)(n-2)(n-3)(n-4)}.
\end{align*} 
\item[8)] Two edges are the same and no pair of edges share any node. 
\begin{align*}
f_{8a} = & (s-r)(s-r-1)(s-2)(n-s-1)(2(s-3)(n-s)+2(t-s)(n-s-2))\\
& +(t-s)(t-s-1)(n-s-2)(n-s-3)((s-r)(s-r-1)+2s(s-1))\\
& +s(s-1)(s-2)(n-s-1)((s-3)(n-s)+2(t-s)(n-s-2)) \\
& +2(s-r)(t-s)(s-1)(s-2)(n-s-1)(n-s-2),\\
f_{8b} = & r(n-s)(n-s-1)(n-r-3)((r-1)(n-r-2) + 2(s-r)(n-r-4)) \\
 & +(s-r)(s-r-1)(n-s)(n-s-1)(n-r-4)(n-r-5)\\
& +(t-s)(t-s-1)(r(r-1)(n-r-2)(n-r-3)\\
& +2r(s-r)(n-r-3)(n-r-4) +(s-r)(s-r-1)(n-r-4)(n-r-5)), \\
P_2 = & \frac{f_{8a}+f_{8b}}{n(n-1)(n-2)(n-3)(n-4)(n-5)}.
\end{align*} 
\item[9)] The four edges share one node, and none of them share the other node (star-shaped). 
$$ P_2 = 0.$$
\item[10)] Linear chain of edges such that two distinct edges share one node with the other two edges and share a node with each other other. 
\begin{align*}
f_{10a} = & (t-s)(s-1)[(s-r)(s-2)(n-s-1)+(s-r)(t-s-1)(n-s-2)\\
& +(n-s-1)(s-r)(s-2)+s(n-s-1)(n-s-2)] \\
f_{10b} = & (s-r)(n-s)(n-s-1)[r(s-2)+(s-r-1)(s-2)]\\
& +(s-r)(t-s)(s-1)[(n-s-1)(n-s-2) + (t-s-1)(n-s-2)]\\
& +(s-r)(t-s)(n-s-1)[r(n-r-3)\\
& +(s-r-1)(n-r-4)+(s-r-1)(s-2)],\\
f_{10c} = & (s-1)(n-s-1)[(s-r)(s-2)(n-s)+(t-s)s(n-s-2)],\\
P_2 = & \frac{f_{10a}+f_{10b}+f_{10c}}{n(n-1)(n-2)(n-3)(n-4)}.
\end{align*} 
\item[11)] All four edges form a box. 
\begin{align*} P_2 = \frac{(s-r)(t-s)(s-1)(n-s-1)}{n(n-1)(n-2)(n-3)}.
\end{align*} 
\item[12)] Three edges form a triangle and one edge connects to one node of the triangle. 
\begin{align*} P_2 = \frac{(s-r)(t-s)(s-1)(n-s-1)}{n(n-1)(n-2)(n-3)}.
\end{align*} 
\item[13)] Three edges share the same node and the fourth edge shares the other node of one of the edges. 
\begin{align*} 
P_2 = & \frac{(s-r)(t-s)(s-1)((s-2)(n-s-1)+(t-s-1)(n-s-2))}{n(n-1)(n-2)(n-3)(n-4)} \\
& +\frac{(s-r)(t-s)(n-s-1)(3r(n-r-3)+2(s-r-1)(n-r-4))}{n(n-1)(n-2)(n-3)(n-4)}\\
& +\frac{(t-s)(r(r-1)(n-s-1)(n-r-2))}{n(n-1)(n-2)(n-3)(n-4)}\\
& +\frac{(s-r)(s-1)(n-s-1)((s-2)(n-s)+(t-s)(n-s-2))}{n(n-1)(n-2)(n-3)(n-4)}.
\end{align*} 
\item[14)] Three edges share the same node and the fourth edge does not share any node with the other edges. 
\begin{align*} 
f_{14a} = & (s-r)(t-s)(s-1)(s-2)((s-3)(n-s-1)+(t-s-1)(n-s-2))\\
& +(t-s)s(s-1)(n-s-2)((n-s-1)(s-2)+(t-s-1)(n-s-3)), \\
f_{14b} = & (s-r)(s-r-1)(s-2)(n-s-1)((t-s)(n-s-2)+(s-3)(n-s))\\
& +(s-r)(s-1)(n-s-1)(n-s-2)((n-s)(s-2)+(t-s)(n-s-3)),\\
P_2 = & \frac{f_{14a} + f_{14b}}{n(n-1)(n-2)(n-3)(n-4)(n-5)}.
\end{align*} 
\item[15)] Three edges form a linear chain and the fourth edge does not share any node with the other edges. 
\begin{align*}
f_{15a} = & (t-s)((n-s-1)(r(r-1)(n-r-2)(n-r-3)\\
& +2r(s-r)(n-r-3)(n-r-4) +(s-r)(s-r-1)(n-r-4)(n-r-5))\\
& +(s-r)(s-1)(s-2)((s-3)(n-s-1) +(t-s-1)(n-s-2)) \\
& +s(s-1)(n-s-2)((n-s-1)(s-2)+(t-s-1)(n-s-3))),\\
f_{15b} = &(s-r)(s-1)(s-2)(n-s-1)((s-3)(n-s)+(t-s)(n-s-2))\\ 
& +(t-s)s(s-1)(n-s-2)((s-2)(n-s-1)+(t-s-1)(n-s-3)), \\
f_{15c} = & (t-s)((s-r)(s-2)((s-r-1)(s-3)(n-s-1)\\
& +(t-s-1)(s-1)(n-s-2)) + (s-1)(n-s-2)((n-s-1)(s-r)(s-2)\\
&  + (t-s-1)s(n-s-3))), \\
f_{15d} = &(s-r)((s-r-1)(s-2)((s-3)(n-s)(n-s-1)\\
& +(t-s)(n-s-1)(n-s-2)) +(t-s)(s-1)((s-2)(n-s-1)(n-s-2)\\
& +(t-s-1)(n-s-2)(n-s-3))+ r(n-s-1)((r-1)(n-s)(n-r-3)\\
&  + 2(s-r-1)(n-s)(n-r-4) + (t-s)(n-s-2)(n-r-4))\\
& + (s-r-1)(n-s-1)(n-r-5)((s-r-2)(n-s) + (t-s)(n-s-2))), \\
f_{15e} = & (s-r)((n-s)(n-s-1)(r(r-1)(n-r-3)+2r(s-r-1)(n-r-4)\\
& +(s-r-1)(s-r-2)(n-r-5)+(s-1)(s-2)(s-3))\\
& + (s-1)(t-s)(n-s-2)(2(s-2)(n-s-1)+(t-s-1)(n-s-3)))\\
& +(t-s)(n-s-2)(n-s-1)(r(r-1)(n-r-3)\\
& +2r(s-r)(n-r-4)+(s-r)(s-r-1)(n-r-5)), \\
P_2 = & \frac{f_{15a} + f_{15b}+f_{15c}+f_{15d}+f_{15e}}{n(n-1)(n-2)(n-3)(n-4)(n-5)}.
\end{align*} 
\item[16)] Three edges form a triangle and the fourth edge does not share any node with the other edges. 
\begin{align*} 
P_2 = & \frac{(t-s)(s-1)((s-r)(s-2)(n-s-1)+(t-s-1)s(n-s-2))}{n(n-1)(n-2)(n-3)(n-4)}\\
& +\frac{(s-r)(n-s-1)((s-r-1)(s-2)(n-s)+(t-s)(s-1)(n-s-2))}{n(n-1)(n-2)(n-3)(n-4)},
\end{align*} 
\item[17)] Two pairs of edges share one node with each other. The pairs of edges do not share any nodes with each other. 
\begin{align*}
P_2 = & \frac{(s-r)(t-s)(s-r-1)(s-2)((s-3)(n-s-1)+(t-s-1)(n-s-2))}{n(n-1)(n-2)(n-3)(n-4)(n-5)} \\
& +\frac{(s-r)(t-s)(s-1)(n-s-2)((n-s-1)(s-2)+(t-s-1)(n-s-3))}{n(n-1)(n-2)(n-3)(n-4)(n-5)} \\
& + \frac{(s-r)(s-1)(s-2)(n-s-1)((s-3)(n-s)+(t-s)(n-s-2))}{n(n-1)(n-2)(n-3)(n-4)(n-5)} \\
& + \frac{s(s-1)(n-s-1)(n-s-2)((s-2)(n-s)+(t-s)(n-s-3))}{n(n-1)(n-2)(n-3)(n-4)(n-5)}. 
\end{align*} 
\item[18)] Two edges share one node with each other. The other edges do not share any nodes with any of the other edges.
\begin{align*} 
f_{18a} = &(s-r)((s-r-1)(s-2)((s-3)(s-4)(n-s)(n-s-1)\\
&+2(t-s)(s-3)(n-s-1)(n-s-2)\\
& +(t-s)(t-s-1)(n-s-2)(n-s-3)) \\
& +(s-1)(n-s-2)((s-2)(s-3)(n-s)(n-s-1)\\
& +2(t-s)(s-2)(n-s-1)(n-s-3)\\
& +(t-s)(t-s-1)(n-s-3)(n-s-4))),\\
f_{18b} = & (s-r)(s-1)(s-2)[(s-3)(s-4)(n-s)(n-s-1)\\
& +2(t-s)(s-3)(n-s-1)(n-s-2)\\
& +(t-s)(t-s-1)(n-s-2)(n-s-3)]\\
& +s(s-1)(n-s-2)((s-2)(s-3)(n-s)(n-s-1)\\
& +2(t-s)(s-2)(n-s-1)(n-s-3)\\
& +(t-s)(t-s-1)(n-s-3)(n-s-4)),\\
f_{18c} = & (s-r)(s-2)((s-r-1)(s-3)(n-s-1)((s-4)(n-s)\\
&+(t-s)(n-s-2))\\
& +(n-s-2)(n-s-1)(s-1)((t-s)(n-s-3)+(s-3)(n-s)))\\
& +(t-s)(s-1)((s-r)(s-2)(n-s-2)((s-3)(n-s-1)\\
& +(t-s-1)(n-s-3))+(n-s-2)(n-s-3)(s(n-s-1)(s-2)\\
& +s(t-s-1)(n-s-4))),\\
f_{18d} = &(s-r)(s-r-1)(s-2)(s-3)(n-s-1)((s-4)(n-s)\\
& +(t-s)(n-s-2)) \\
& +2(s-r)(s-1)(s-2)(n-s-1)(n-s-2)((n-s)(s-3)\\
& +2(t-s)(n-s-3)) \\
& +s(s-1)(n-s-1)(n-s-2)(n-s-3)((n-s)(s-2)\\
& +(t-s)(n-s-4)),\\
f_{18e} = &(t-s)((s-r)(s-r-1)(s-2)(s-3)((s-4)(n-s-1)\\
& +(t-s-1)(n-s-2)) \\
& +2(s-r)(s-1)(s-2)(n-s-2)((n-s-1)(s-3)\\
& +(t-s-1)(n-s-3)) +s(s-1)(n-s-1)(n-s-3)((n-s-1)(s-2)\\
& +(t-s-1)(n-s-4))),\\
P_2 = & \frac{f_{18a} + f_{18b}+f_{18c}+f_{18d}+f_{18e}}{n(n-1)(n-2)(n-3)(n-4)(n-5)(n-6)}.
\end{align*} 
\item[19)] None of the four edges share any node. 
\begin{align*} 
f_{19a} = & (s-r)(s-r-1)(s-2)(s-3)((s-4)(s-5)(n-s)(n-s-1)\\
& +2(t-s)(s-4)(n-s-1)(n-s-2)\\
& +(t-s)(t-s-1)(n-s-2)(n-s-3)), \\
f_{19b} = & 2(s-r)(s-1)(s-2)(n-s-2)((n-s)(s-3)(s-4)(n-s-1)\\
& +2(n-s-1)(t-s)(s-3)(n-s-3)\\
& +(t-s)(t-s-1)(n-s-3)(n-s-4)),\\
f_{19c} = & s(s-1)(n-s-2)(n-s-3)((n-s)(n-s-1)(s-2)(s-3)\\
& +2(n-s-1)(t-s)(s-2)(n-s-4)\\
& +(t-s)(t-s-1)(n-s-4)(n-s-5)), \\
P_2 = & \frac{f_{19a}+f_{19b} + f_{19c}}{n(n-1)(n-2)(n-3)(n-4)(n-5)(n-6)(n-7)}.\\
\end{align*} 
\end{enumerate} 

For the remaining expressions, similar derivations using combinatorial analysis can be obtained. 
\subsection{Expression for $Z_w$}
The similarity graph $G$ can be a generic graph constructed from a similarity measure, such as the Euclidean distance. Without loss of generality, $|G| = O(kn)$ with $k = O(n^\alpha), 0 \le \alpha <1$. We assume that $u,v,w = O(1)$. To establish (\ref{eq:condZw}), we focus on the leading terms on the left-hand side of the inequality. After extensive simplification, the leading term for the denominator of $\boldE ((Z^n_w(v) - Z^n_w(u))^{2} (Z^n_w(w) - Z^n_w(v))^{2})$  is
\begin{align} 
\text{den}_{Z_w} \triangleq v^2w^2(k n^2 - \sum_{i=1}^n|G_i|^2)^2(1 - u)^2(1 - v)^2. \label{ref:den_Zw}
\end{align} 

The leading term for the numerator is: %(with respect to $n$)  

\begin{align*} 
& \text{num}_{Z_w} \triangleq (w-v)(v-u)\Big(k^2n^4C_{w,1} + x_{14}C_{w,2} + C_{w,3}\sum_{i=1}^n |G_i|^4 + nC_{w,4} \sum_{i=1}^n |G_i|^3 \\
& + C_{w,5} \sum_i \sum_{j \in G_i; j \neq i} (|G_i|-1)^2(|G_j|-1) + k n^2 C_{w,6} \sum_{i=1}^n |G_i|^2 + n^2 C_{w,7}\sum_{i=1}^n |G_i|^2 \\
& + knC_{w,8}\sum_{i,j \in G} (|G_i|-1)(|G_j|-1) + nC_{w,9}\sum_{i,j \in G} (|G_i|-1)(|G_j|-1) \\
& + nx_7C_{w,10} + n^x_8 C_{w,11} + nx_9 C_{w,12}\Big)
\end{align*} 
with 
\begin{align*} 
C_{w,1} = & 4vw(1-v)(1-u)  + 2(v - u)(w - v),\\
C_{w,2} = & 8vw(v-u)(1 - u)(1-v),\\
C_{w,3} = & -2(v - u)(w - v) + 2v(1 - u)(1 + v) + vw(5u - 7)(1 - v), \\
C_{w,4} = & 8v(w-v) - 8w +  2v(2 + 9w)(1 - u)(1 - v), \\
C_{w,5} = & 8(w - uv) + (48 - 56v)(w - v) + 16(3v^2 + w)(1 - u)  - 4vw(49 - 37u)(1 - v), \\
C_{w,6} = & -4(v-u)(w-v)  - 8vw(1-v)(1-u),\\
C_{w,7} =& 2(w - uv) + 2(1-2v)(w - v)  + vw(9u - 11)(1 - v) + 2v^2(1-u), \\
C_{w,8} = & 16(v-u)(w-v)  + 32vw(1-v)(1-u),\\
C_{w,9} = & 2(28v - 23)(w - v) - 2(23v^2 + 9w)(1 - u) - 2vw(72u - 95)(1 - v)+ 10(uv-w), \\
C_{w,10} = & -8vw(1 - u)(1 - v)-4(w-v)(v - u),\\
C_{w,11} = & 4vw(1 - u)(1 - v) +2(w-v)(v - u),\\
C_{w,12} = & 8v(5v(1-v) - (1 - u)(12v^2 - 7v + 2))\\
&  - (w-v)(24(1-u)+8(1 - v)(12uv - 17v + 4)).
\end{align*} 
Since $(w-v)(v-u) < (w-u)^2$ for $u<v<w$, the expression $\text{num}_{Z_w}/\text{den}_{Z_w}$ can be bounded by $C(w-u)^2$ as long as the ratio of graph configurations in the numerator and denominator can be bounded asymptotically by $O(1)$. Specifically, since $u,v,w = O(1)$, the terms $C_{w,1}, \hdots, C_{w,12}$ can be bounded asymptotically by a constant. The remaining terms in the numerator involve configurations of the graph: $k$, $n$, $\sum_{i=1}^n |G_i|^4$, $\sum_{i=1}^n |G_i|^3 $, $\sum_{i=1}^n |G_i|^2$, $\sum_i \sum_{j \in G_i; j \neq i} (|G_i|-1)^2(|G_j|-1)$, $\sum_{i,j \in G} (|G_i|-1)(|G_j|-1)$, $x_7$, $x_8$, and $x_9$. If the ratio of each of these terms with the denominator's $(k n^2 - \sum_{i=1}^n|G_i|^2)^2$ is bounded by $O(1)$, then the entire expression can be asymptotically bounded by a constant $C_w$ times $(w-u)^2$.

In the following, we assume that $\sum_{i=1}^n|G_i|^2 = o(kn^2)$ and we check each configuration (in their order of appearance).% is of smaller order than $kn^2$. 

Clearly  $\frac{k^2n^4}{(k n^2 - \sum_{i=1}^n|G_i|^2)^2} \precsim O(1)$.  

For $x_{14}$, we have
\begin{align*}
x_{14} & = \sum_i \sum_{j \neq i} (|G_i \setminus \{j \in G_i\} |) (|G_i \setminus \{j \in G_i\}|-1)(|G_j \setminus \{i \in G_j\}|)(|G_j \setminus \{i \in G_j\}|-1) \\
& < \sum_i \sum_{j \neq i} |G_i|^2 |G_j|^2 \\
& = |G_1|^2 \sum_{j \neq 1} |G_j|^2  + |G_2|^2 \sum_{j \neq 2} |G_j|^2 + \hdots + |G_n|^2 \sum_{j \neq n} |G_n|^2 \\
& =  |G_1|^2 (\sum_{i=1}^n |G_i|^2 - |G_1|^2)  + |G_2|^2 (\sum_{i=1}^n |G_i|^2 - |G_2|^2)  + \hdots + |G_n|^2 (\sum_{i=1}^n |G_i|^2 - |G_n|^2)  \\
& = (\sum_{i=1}^n |G_i|^2)^2 - \sum_{i=1}^n |G_i|^4 >0.
\end{align*} 
Then $x_{14} < (\sum_{i=1}^n |G_i|^2)^2$ and $\frac{x_{14}}{(k n^2 - \sum_{i=1}^n|G_i|^2)^2} \precsim O(1)$. Following similar arguments, since $\sum_{i=1}^n |G_i|^2 = o(kn^2)$, we have  $\frac{\sum_{i=1}^n |G_i|^4}{(k n^2 - \sum_{i=1}^n|G_i|^2)^2} \precsim O(1) $ and $\frac{n\sum_{i=1}^n |G_i|^3}{(k n^2 - \sum_{i=1}^n|G_i|^2)^2} \precsim O(1)$. 

For $\sum_i \sum_{j \in G_i; j \neq i} (|G_i|-1)^2(|G_j|-1)$, we have
 $$\sum_i \sum_{j \in G_i; j \neq i} (|G_i|-1)^2(|G_j|-1) < 2|G|\sum_{i=1}^n |G_i|^3.$$
 Since the the largest $|G_i|$ can be is $n-1$ (every other observation connects to node $\textbf{y}_i$), it follows that $2|G|\sum_{i=1}^n |G_i|^3 \precsim 2|G|kn^3 \asymp k^2n^4 $
and $\frac{k^2n^4}{(k n^2 - \sum_{i=1}^n|G_i|^2)^2} \precsim O(1)$. 

Similarly, since $k n^2\sum_{i=1}^n |G_i|^2 \precsim k^2n^4 $, we have $\frac{kn^2 \sum_{i=1}^n |G_i|^2 }{(k n^2 - \sum_{i=1}^n|G_i|^2)^2} \precsim O(1).$

We have $\sum_{(i,j) \in G} (|G_i|-1)(|G_j|-1) < \sum_{i=1}^n |G_i| (|G| - |G_i|) = |G| \sum_{i=1}^n |G_i| - \sum_{i=1}^n |G_i|^2 < 2|G|^2 \asymp 2k^2n^2$, and so  $\frac{\sum_{i,j \in G} (|G_i|-1)(|G_j|-1)}{(k n^2 - \sum_{i=1}^n|G_i|^2)^2} \precsim O(1).$

Finally, since 
\begin{align*} 
nx_7 & < n \sum_{i=1} |G_i| (|G|-|G_i|) < n2|G|^2 \asymp k^2n^3 ,\\
x_8 &= \sum_{(i,j),(j,l), i\neq l} |\{ l: (i,l), (l,m) \in G\} \precsim kn^3, \\
x_9 & = \sum_{(i,j)} \sum_{l:(i,l),(j,l) \in G} (|G_l| - 2) \precsim k^2n^4 ,
\end{align*}
it follows that the ratio of the these configurations with $(k n^2 - \sum_{i=1}^n|G_i|^2)^2$  are bounded asymptotically by $O(1)$. 

\subsection{Expression for $Z_\text{diff}$}

We adopt a similar approach for $Z_\text{diff}$: we study the analytical expression for $\boldE\left((Z^n_\text{diff}(v) - Z^n_\text{diff}(u))^{2} (Z^n_\text{diff}(w) - Z^n_\text{diff}(v))^{2}\right)$. This expression can be written as the combination of terms involving $u,v,$ and $w$ and terms involving configurations from the graph. We first show that the expressions involving $u,v,$ and $w$ can be bounded by $C(w-u)^2$ or $C(w-u)$. We then show that the graph-configurations are bounded asymptotically by $O(1)$ or $O(1/n)$.  It follows then that the entire expression can be bounded by a constant $C_\text{diff}$ times $(w-u)^2$. 

Let $e_v = v(1-v)$, $e_w = w(1-w)$, and $e_u = u(1-u)$. 
The leading term for the denominator of $\boldE\left((Z^n_\text{diff}(v) - Z^n_\text{diff}(u))^{2} (Z^n_\text{diff}(w) - Z^n_\text{diff}(v))^{2}\right)$ is: 
\begin{align*} 
\text{den}_{Z_\text{diff}} & = (nV_G)^2w(1-u)e_u e_v^3 e_w
\end{align*} 
with $V_G = \sum_i |G_i|^2 - 4|G|^2/n$. 

For the numerator of $\boldE\left((Z^n_\text{diff}(v) - Z^n_\text{diff}(u))^{2} (Z^n_\text{diff}(w) - Z^n_\text{diff}(v))^{2}\right)$, we group the leading terms by their graph configurations. The numerator can be expressed as
\begin{align*} & K_1(u,v,w) \times k^4 n^2 + K_2(u,v,w) \times k^2n (\sum_{i=1}^n |G_i|^2) +  K_3(u,v,w) \times \sum_{i=1}^n |G_i|^4 \\
& + K_4(u,v,w) \times k \sum_{i=1}^n |G_i|^3 +  K_5(u,v,w)\times x_{14} \\
& +  K_6(u,v,w) \times \sum_i \sum_{j \in G_i; j \neq i} (|G_i|-1)^2(|G_j|-1).
\end{align*}
We first show that the coefficients $K_1(u,v,w), K_2(u,v,w), K_3(u,v,w), K_4(u,v,w)$, and $K_5(u,v,w)$ can be bounded by $C(w-u)^2$ or $C(w-u)$. 
 
\begin{enumerate} 
\item $K_1(u,v,w)$: The leading coefficient for $k^4 n^2$ can be expanded as 
\begin{align*} 
K_1(u,v,w) = & C_{d,1} (w-v)^2 + C_{d,2} (v-u)(w-v) + C_{d,3}\sqrt{e_u}(\sqrt{e_u}-\sqrt{e_v})(w-v) \\
& + C_{d,4}\sqrt{e_v}(\sqrt{e_v} - \sqrt{e_w})\Big(\sqrt{u(1-v)}(\sqrt{v(1-u)} - \sqrt{u(1-v)})\\
&  -2\sqrt{v(1-u)}(\sqrt{v(1-u)}-\sqrt{u(1-v)})\Big)
%& 32 \times\Big\{ 2v(1 - v)\sqrt{e_u}\left[\sqrt{u(1-v)}(\sqrt{u(1-v)} -\sqrt{v(1-u)}) \right. \\
%& \left. -2\sqrt{v(1-u)}\left(\sqrt{u(1-v)}-\sqrt{v(1-u)}\right)\right]\blue{(w-v)^2} \\
%& + \sqrt{e_u}(12u - 8)\blue{(v-u)^4(w-v)}\\
%& + 6\sqrt{e_u}(8u^2 - 9u + 2 + \sqrt{e_u}(2\sqrt{e_v} - \sqrt{e_w}))\blue{(v-u)^3(w-v)} \\
%& + [2\sqrt{e_u}(36u^3 - 57u^2 + 24u - 2)-18u\sqrt{e_v}(1-u)(1-2u) + 2u\sqrt{e_w}(1-u)(5-9u) \\
%& + 2\sqrt{e_u}\sqrt{e_v}\sqrt{e_w}(3u - 2)]\blue{(v-u)^2(w-v)} \\
%& + [-2u\sqrt{e_u}(24u^2 - 25u + 5)(1-u) + 6u\sqrt{e_v}(1-u)(6u^2 - 6u + 1)\\
%&  - 2u\sqrt{e_w}(1-u)(9u^2 - 10u + 2) + \sqrt{e_u}\sqrt{e_v}\sqrt{e_w}(12u^2 - 14u + 3)]\blue{(v-u)(w-v)} \\
%& +2\blue{u(1 - u)(\sqrt{e_u}-\sqrt{e_v})}(\sqrt{e_u}\sqrt{e_w}(3u - 2) -3u(1-u)(1-2u) )\blue{(w-v)} \\
%& + 2v(1-v)\sqrt{e_u}\blue{\sqrt{e_v}(\sqrt{e_v} - \sqrt{e_w})}[\blue{\sqrt{u(1-v)}(\sqrt{v(1-u)} - \sqrt{u(1-v)})} \\
%& -2\blue{\sqrt{v(1-u)}(\sqrt{v(1-u)}-\sqrt{u(1-v)})}] \Big\}.
\end{align*} 
with
\begin{align*} 
C_{d,1} = & 64v(1 - v)\sqrt{e_u}(u(1-v)+2v(1-u)- 3\sqrt{e_ue_v}),\\
C_{d,2} = & 32 \sqrt{e_u}(12u - 8)(v-u)^3 + 192\sqrt{e_u}(8u^2 - 9u + 2 + \sqrt{e_u}(2\sqrt{e_v} - \sqrt{e_w}))(v-u)^2 \\
& + 32(2\sqrt{e_u}(36u^3 - 57u^2 + 24u - 2)-18u\sqrt{e_v}(1-u)(1-2u) \\
& + 2u\sqrt{e_w}(1-u)(5-9u) + 2\sqrt{e_u}\sqrt{e_v}\sqrt{e_w}(3u - 2))(v-u) \\
& -64u\sqrt{e_u}(24u^2 - 25u + 5)(1-u) + 192u\sqrt{e_v}(1-u)(6u^2 - 6u + 1) \\
&- 64u\sqrt{e_w}(1-u)(9u^2 - 10u + 2)  + 32\sqrt{e_u}\sqrt{e_v}\sqrt{e_w}(12u^2 - 14u + 3),\\
C_{d,3} = & 64e_u(\sqrt{e_w}(3u - 2) -3\sqrt{e_u}(1-2u)),\\
C_{d,4} = & 64e_v\sqrt{e_u}.
\end{align*} 
It is clear that $C_{d,1} (w-v)^2 + C_{d,2} (v-u)(w-v) \le C(w-u)^2$ since
$C_{d,1} (w-v)^2 + C_{d,2} (v-u)(w-v) \le (C_{d,1} + C_{d,2}) (w-u)^2$ and $C$ can be chosen to be large enough such that $C_{d,1} + C_{d,2} \le C$. 
In the following we focus on the next two terms. For the third term, we need to show that $\sqrt{e_u}(\sqrt{e_u}-\sqrt{e_v}) \le (v-u)$. Let $\delta = v-u$ and define
\begin{align*} 
 g(\delta) & = \sqrt{e_u}(\sqrt{e_u}-\sqrt{e_v}) \\
& = \sqrt{u(1-u)}\left(\sqrt{u(1-u)} - \sqrt{(u+\delta)(1-u-\delta)} \right)
\end{align*}
which is continuous everywhere on $0 \le \delta \le 1-u$. 

If $g(\delta)$ is convex for $0 \le \delta \le 1-u$, it follows that $g(\delta) \le \delta$. Since $g(0) = 0$ and $g(1-u) = u(1-u) \le 1(-u)$, what remains is to check its second derivative is non-negative:
\begin{align*}
g'(\delta)  & = \frac{-(1-2u-2\delta)\sqrt{u(1-u)}}{2\sqrt{(u+\delta)(1-u-\delta)}},\\
g''(\delta) & = \frac{\sqrt{u(1-u)}}{2} \left(\frac{2}{\sqrt{(u+\delta)(1-u-\delta)}}  + \frac{(1-2u-2\delta)^2}{2\sqrt{(u+\delta)(1-u-\delta)^3}}\right) >0.
\end{align*} 
Since we have established that  $g(\delta) = \sqrt{e_u}(\sqrt{e_u}-\sqrt{e_v})$ is convex, it follows that $\sqrt{e_u}(\sqrt{e_u}-\sqrt{e_v}) \le (v-u)$ and $\sqrt{e_v}(\sqrt{e_v} - \sqrt{e_w}) \le (w-v)$. Moreover, the minimum of $g(\delta)$ is achieved when $\delta = 0.5 - u$ and $-g(0.5-u) = \sqrt{e_u}(\frac{1}{2}-\sqrt{e_u}) \le \frac{1}{2} - u$, for $u<\frac{1}{2}$. Therefore $|\sqrt{e_u}(\sqrt{e_u} - \sqrt{e_v})| \le (v-u)$. 

 Following a similar argument, we can establish that
$\sqrt{u(1-v)}(\sqrt{v(1-u)} - \sqrt{u(1-v)}) \le (v-u)$. Let $h(\delta ) = \sqrt{(u+\delta)(1-u)}\left(\sqrt{(u+\delta)(1-u)} - \sqrt{u(1-u-\delta)} \right)$. We have $h(0) = 0$ and $h(1-u) = 1-u$. 
Its first and second derivatives are
\begin{align*} 
h'(\delta) &= (1-u) - \frac{(1-2u-2\delta)\sqrt{u(1-u)}}{2\sqrt{(1-u-\delta)(u+\delta)}},\\
h''(\delta) & = \frac{1}{2} \sqrt{u(1-u)} \left( \frac{2}{\sqrt{(1-u-\delta)(u+\delta)}} + \frac{(1-2u-2\delta)^2}{2\sqrt{(1-u-\delta)(u+\delta)}^3} \right) >0,
\end{align*} 
and therefore $\sqrt{v(1-u)}(\sqrt{v(1-u)} - \sqrt{u(1-v)}) \le (v-u)$. 
Since $\sqrt{u(1-v)} < \sqrt{v(1-u)}$, it follows that $\sqrt{u(1-v)}(\sqrt{v(1-u)} - \sqrt{u(1-v)}) \le (v-u)$. Note that $\sqrt{v(1-u)} - \sqrt{u(1-v)} >0$. 

Therefore, $K_1(u,v,w) \le C(w-u)^2$ for some constant $C$. 
%Then as long as $k^4n^2/(nV_G)^2$ is at most of $O(1)$, the expression can be bounded by $C(w-u)^2$. 
\item $K_2(u,v,w)$: The leading coefficient for $k^2n(\sum_{i=1}^n |G_i|^2) $ is
\begin{align*}
K_2(u,v,w) = & \sqrt{e_u}\sqrt{e_v}\sqrt{e_w}\Big\{C_{d,5}(v-u)^2 + C_{d,6}(w-u)(v-u) + C_{d,7}(\sqrt{e_u}-\sqrt{e_v})(v-u) \\
& +  C_{d,8}(\sqrt{e_u} - \sqrt{e_w})(v-u) + C_{d,9}\sqrt{e_u}(\sqrt{e_u}-\sqrt{e_v})^2 \\
& + C_{d,10} (\sqrt{e_u}-\sqrt{e_v})(\sqrt{u(1-w)}(\sqrt{w(1-u)} - \sqrt{u(1-w)}) \\
& - 2\sqrt{w(1-u)}(\sqrt{w(1-u)}-\sqrt{u(1-w)}))\Big\}
\end{align*} 
with
\begin{align*} 
C_{d,5} = & 16(-6uw + 4w + 2u - 1)\\
C_{d,6} = & 16(-12u^2 + 14u - 6\sqrt{e_u}\sqrt{e_v} - 3)\\
C_{d,7} = & 16(-2\sqrt{e_w}(2 - 3u) - 2\sqrt{e_u}(1-3u))\\
C_{d,8} = & 32(3u - 2)\sqrt{e_u} \\
C_{d,9} = & -96 \sqrt{e_w} \\
C_{d,10} = & 32\sqrt{e_u}
\end{align*} 
%\begin{align*}
%& \sqrt{e_u}\sqrt{e_v}\sqrt{e_w}\Big\{16(-6uw + 4w + 2u - 1)\blue{(v-u)^2} \\
%&  +16(-12u^2 + 14u - 6\sqrt{e_u}\sqrt{e_v} - 3)\blue{(w-u)(v-u)} \\
%& +16(-2\blue{\sqrt{e_w}}(2 - 3u) - 2\blue{\sqrt{e_u}}(1-3u))\blue{(\sqrt{e_u}-\sqrt{e_v})(v-u)} \\
%& + 32(3u - 2)\blue{\sqrt{e_u}(\sqrt{e_u} - \sqrt{e_w})(v-u)} \\
%& + - 96\blue{\sqrt{e_u}\sqrt{e_w}(\sqrt{e_u}-\sqrt{e_v})^2} \\
%& + 32\blue{\sqrt{e_u}(\sqrt{e_u}-\sqrt{e_v})}[\blue{\sqrt{u(1-w)}(\sqrt{w(1-u)} - \sqrt{u(1-w)})} \\
%& - 2\blue{\sqrt{w(1-u)}(\sqrt{w(1-u)}-\sqrt{u(1-w)})}]\Big\}.
%\end{align*} 

Since $ u < v < w$, we have $C_{d,5}(v_u)^2 + C_{d,6}(w-u)(v-u) \le C(w-u)^2$ for some constant $C$. In order to show that remaining terms can also be bounded by $C(w-u)^2$, we follow that same argument detailed above for $K_1(u,v,w)$. Observe that $|\sqrt{e_u}(\sqrt{e_u} - \sqrt{e_v})| \le (v-u)$ and $|\sqrt{e_u}(\sqrt{e_u} - \sqrt{e_w})| \le (w-u)$. It follows that terms with $C_{d,7}, C_{d,8},$ and $C_{d,9}$ of $K_2(u,v,w)$ can be by bounded by $C(w-u)^2$ as well. 

Finally, for the last term in $K_2(u,v,w)$, we see that $\sqrt{u(1-w)}(\sqrt{w(1-u)} - \sqrt{u(1-w)}) \le (w-u)$ and $\sqrt{w(1-u)}(\sqrt{w(1-u)} - \sqrt{u(1-w)}) \le (w-u)$. 

It follows that $K_2(u,v,w) \le C(w-u)^2$ for some constant $C$. 
%Then as long as $k^2n(\sum_{i=1}^n |G_i|^2)/(nV_G)^2$ is at most of $O(1)$, the expression can be bounded by $C(w-u)^2$. 
\item $K_3(u,v,w)$: The leading coefficient for $\sum_{i=1}^n |G_i|^4$ is
\begin{align*}
K_3(u,v,w) = & C_{d,11}(w-v)^2 + C_{d,12} (v-u)(w-v) + C_{d,13} \sqrt{e_u}(\sqrt{e_u}-\sqrt{e_v})(w-v)\\
& + C_{d,14} \sqrt{e_v}(\sqrt{e_v}-\sqrt{e_u})(w-v) \\
& + C_{d,15} \sqrt{e_v}(\sqrt{e_v}-\sqrt{e_w})\sqrt{e_u}(\sqrt{e_u}-\sqrt{e_v}) + C_{d,16} \sqrt{e_v}(\sqrt{e_v}-\sqrt{e_w})(v-u)
\end{align*} 
with 
\begin{align*}
C_{d,11} = & -2v\sqrt{e_u}(3u + 5v - 8uv + 6uv^2 - 4v^2 - 2) - 4e_u\sqrt{e_v}(3v^2 - 3v + 1),\\
C_{d,12} = &  8\sqrt{e_u}(2 - 3u)(v-u)^3(w-v), \\
& - (4\sqrt{e_u}(24u^2 - 27u + 7) + 24e_u\sqrt{e_v} + 12e_u\sqrt{e_w})(v-u)^2(w-v), \\
& +  (-2\sqrt{e_u}(72u^3 - 114u^2 + 56u - 9)+ 36u\sqrt{e_v}(2u^2 - 3u + 1), \\
& -4u\sqrt{e_w}(9u^2 - 14u + 5) -4\sqrt{e_u}\sqrt{e_v}\sqrt{e_w}(3u - 2))(v-u)(w-v),\\
& -2\sqrt{e_u}(48u^4 - 98u^3 + 70u^2 - 21u + 2) + 4u\sqrt{e_v}(18u^3 - 36u^2 + 23u - 5),\\
&   -2u\sqrt{e_w}(18u^3 - 38u^2 + 25u - 5)-\sqrt{e_u}\sqrt{e_v}\sqrt{e_w}(24u^2 - 28u + 7), \\
&  -2\sqrt{e_w}(1-u-v)(6u^3 - 10u^2 + 5u - 1),\\
C_{d,13} = &  -4\sqrt{e_u}(14u^2 - (6u-1)(u^2+u+1)),  \\
C_{d,14} = &  2\sqrt{e_w}(6u^3 - 10u^2 + 5u - 1),\\
C_{d,15} = &  4\sqrt{e_v}(3v^2-3v+1),\\
C_{d,16} = & 2v\sqrt{e_u}(6v^2 - 8v + 3).
\end{align*} 
%\begin{align*} 
%&[-2v\sqrt{e_u}(3u + 5v - 8uv + 6uv^2 - 4v^2 - 2) - 4e_u\sqrt{e_v}(3v^2 - 3v + 1)]\blue{(w-v)^2} \\
%& + 8\sqrt{e_u}(2 - 3u)\blue{(v-u)^4(w-v)} \\
%& + (-4\sqrt{e_u}(24u^2 - 27u + 7) - 24e_u\sqrt{e_v} + 12e_u\sqrt{e_w})\blue{(v-u)^3(w-v)} \\
%& + (-2\sqrt{e_u}(72u^3 - 114u^2 + 56u - 9)+ 36u\sqrt{e_v}(2u^2 - 3u + 1) \\
%& -4u\sqrt{e_w}(9u^2 - 14u + 5) -4\sqrt{e_u}\sqrt{e_v}\sqrt{e_w}(3u - 2))\blue{(v-u)^2(w-v)}\\
%& + (-2\sqrt{e_u}(48u^4 - 98u^3 + 70u^2 - 21u + 2) + 4u\sqrt{e_v}(18u^3 - 36u^2 + 23u - 5)\\
%&  -2u\sqrt{e_w}(18u^3 - 38u^2 + 25u - 5)-\sqrt{e_u}\sqrt{e_v}\sqrt{e_w}(24u^2 - 28u + 7))\blue{(v-u)(w-v)} \\
%& + ( -4\blue{e_u}(14u^2 - (6u-1)(u^2+u+1)) -2\blue{\sqrt{e_v}}\sqrt{e_w}(6u^3 - 10u^2 + 5u - 1))\blue{(\sqrt{e_u}-\sqrt{e_v})(w-v)} \\
%& -2\sqrt{e_w}(1-u-v)(6u^3 - 10u^2 + 5u - 1)\blue{(v-u)(w-v)} \\
%& +2\blue{\sqrt{e_u}}\blue{(\sqrt{e_v} - \sqrt{e_w})}(2\blue{e_v}(3v^2-3v+1)\blue{(\sqrt{e_u}-\sqrt{e_v})} - v\sqrt{e_v}\blue{(v - u)}(6v^2 - 8v + 3)).
%\end{align*} 
Again, the first two terms involving $C_{d,11}$ and $C_{d,12}$ can be bounded by $C(w-u)^2$. Repeating the convexity argument, $|\sqrt{e_v}(\sqrt{e_v} - \sqrt{e_u})| \le (v-u)$, which allow us to bound the remaining terms by $C(w-u)^2$ as well.   Therefore, the entire expression $K_3(u,v,w)$ can also be bounded by $C(w-u)^2$. \\
\item $K_4(u,v,w)$: The leading coefficient for $k\sum_{i=1}^n |G_i|^3$ is
\begin{align*} 
& C_{d,17}(w-v)^2+ C_{d,18}(w-v)(v-u) + C_{d,19}\sqrt{e_u}(\sqrt{e_u}-\sqrt{e_v})(w-v) + C_{d,20}(v-u)^2 \\
& + C_{d,21}\sqrt{e_u}(\sqrt{e_v} - \sqrt{e_w})(v-u) + C_{d,22}e_u(\sqrt{e_u}-\sqrt{e_v})(\sqrt{e_v}-\sqrt{e_w})
\end{align*} 
with 
\begin{align*} 
C_{d,17} = & \sqrt{e_u}(16uv(12v^2 - 16v + 5) -16v(8v^2 - 9v + 2))+ 32e_u\sqrt{e_v}(6v^2 - 6v + 1),\\
C_{d,18} = & (128\sqrt{e_u}(3u - 2))(v-u)^3 \\
& + (32\sqrt{e_u}(48u^2 - 54u + 13) + 384e_u\sqrt{e_v} - 192e_u\sqrt{e_w} )(v-u)^2\\
&  + (16\sqrt{e_u}(144u^3 - 228u^2 + 104u - 13)  -576\sqrt{e_v}u(2u^2 - 3u + 1) \\
& + 64\sqrt{e_w}u(9u^2 - 14u + 5) + 64\sqrt{e_u}\sqrt{e_v}\sqrt{e_w}(3u - 2))(v-u)\\
&  + 16\sqrt{e_u}(96u^4 - 196u^3 + 130u^2 - 31u + 2), \\
C_{d,19} = &  8(\sqrt{e_u}(1 - 2u)(4(1-6u(1-u)))  + 2\sqrt{e_w}(1-u)(1-8u+3u^2)),\\
C_{d,20} = & (\sqrt{e_v}-\sqrt{e_w})(-192e_u(v-u)^2 \\
& + 64(6u - 2\sqrt{e_u}\sqrt{e_v} - 18u^2 + 12u^3 + 3\sqrt{e_u}\sqrt{e_v}u)(v-u) \\
&  - 32e_u(7 - 36u(1-u)) + 16\sqrt{e_u}\sqrt{e_v}(9 - 40u + 36u^2)), \\
C_{d,21} = & 32\sqrt{e_u}(1-2u)(1- 12u + 12u^2) \\
& + 16\sqrt{e_v}(u(36u^2 - 56u + 23) - 2 + 5u  - 14u^2 + 9u^3) \\
C_{d,22} = & 32\sqrt{eu}(1 - 6u + 6u^2) 
\end{align*} 
The first two terms involving $C_{d,17}$, $C_{d,18}$, and $C_{d,20}$ can be bounded by $C(w-u)^2$. Repeating a combination of the convexity arguments from above, the remaining terms can also be bounded by $C(w-u)^2$. It follows that $K_4(u,v,w) \le C(w-u)^2$.   
\item $K_5(u,v,w)$: The leading coefficient for $x_{14}$ is
\begin{align*}
K_5(u,v,w) = & C_{d,23} (w-v)^2 + C_{d,24} \sqrt{e_u}(\sqrt{e_u}-\sqrt{e_v})(w-v) + C_{d,25}(v-u)(w-v) \\
& + C_{d,26}\sqrt{e_v}(\sqrt{e_v} - \sqrt{e_w})(v-u) + C_{d,27}\sqrt{e_v}(\sqrt{e_v} - \sqrt{e_w})\sqrt{e_u}(\sqrt{e_u} - \sqrt{e_v}).
\end{align*} 
with
\begin{align*}
C_{d,23} = & 4v(1 - v)((u + 2v - 3uv)\sqrt{e_u} - 3\sqrt{e_v}u(1-u))\\ 
C_{d,24} = & -2\sqrt{e_v}((6u\sqrt{e_w} - 6\sqrt{e_u} + 12\sqrt{e_u}v)(\sqrt{e_u} - \sqrt{e_v})\\
&  -2(1 - 2v)(u + 2v - 3uv)+ 12u(1-u)(1 - 2v) + 2\sqrt{e_u}\sqrt{e_w}(-6u + 3v - 2))\\
C_{d,25} = & -2\sqrt{e_u}\sqrt{e_v}((4\sqrt{e_w} - 8\sqrt{e_u} + 12\sqrt{e_u}u)(v-u)\\
& +2\sqrt{e_u}(2u - 1)(3u - 2) + \sqrt{e_w}(8u - 3)\\
&  - 6\sqrt{e_u}^2\sqrt{e_w})\\
C_{d,26} = & - 4\sqrt{e_v}((- 3u^2 + 3u)(v-u) + 2\sqrt{e_u}\sqrt{e_v} - 3u + 9u^2 - 6u^3 - 3u\sqrt{e_u}\sqrt{e_v})\\
C_{d,27} = & 12\sqrt{e_v}u(1-u).
\end{align*} 
%\begin{align*} 
%& \Big[4v(1 - v)((u + 2v - 3uv)\sqrt{e_u} - 3\sqrt{e_v}u(1-u))\Big]\blue{(w-v)^2}\\
%&  -2\blue{\sqrt{e_u}\sqrt{e_v}}[(6\sqrt{e_w}u - 6\sqrt{e_u} + 12\sqrt{e_u}v)\blue{(\sqrt{e_u}-\sqrt{e_v})^2} \\
%& + \blue{(\sqrt{e_u}-\sqrt{e_v})} \left(-2(1 - 2v)(u + 2v - 3uv)+ 12u(1-u)(1 - 2v) + 2\sqrt{e_u}\sqrt{e_w}(-6u + 3v - 2)\right)\\
%& +  (4\sqrt{e_w} - 8\sqrt{e_u} + 12\sqrt{e_u}u)\blue{(v-u)^2} + (2\sqrt{e_u}(2u - 1)(3u - 2) + \sqrt{e_w}(8u - 3))\blue{(v-u)} \\
%& - 6\sqrt{e_u}^2\sqrt{e_w}(v - u)]\blue{(w-v)} \\
%&  - 4\blue{e_v(\sqrt{e_v} - \sqrt{e_w})}\Big[(- 3u^2 + 3u)\blue{(v-u)^2} + (2\sqrt{e_u}\sqrt{e_v} - 3u + 9u^2 - 6u^3 - 3\sqrt{e_u}\sqrt{e_v}u)\blue{(v-u)}\\
%&  - 3u(1-u)\blue{\sqrt{e_u}(\sqrt{e_u} - \sqrt{e_v})}\Big]
%\end{align*} 
and utilizing the arguments above, this term is also bounded by $C(w-u)^2$. 
\item $K_6(u,v,w)$: The leading coefficient for $\sum_i \sum_{j \in G_i; j \neq i} (|G_i|-1)^2(|G_j|-1)$ is
\begin{align*} 
K_6(u,v,w) = C_{d,28} (v-u) +C_{d,29}\sqrt{e_u}( \sqrt{e_u}-\sqrt{e_v})
\end{align*} 
with
\begin{align*}
C_{d,28} = & 16\sqrt{e_u}\sqrt{e_w}(\sqrt{e_u} + 2\sqrt{e_w} - 3\sqrt{e_w}u - 3\sqrt{e_u}w)(v-u)^2 \\
&  - 8\sqrt{e_u}\sqrt{e_w}(-2\sqrt{e_u}(3u + 5w - 9uw - 1)+ \sqrt{e_v}(2u + 4w - 6uw - 1)\\
& + 2\sqrt{e_w}(9u^2 - 10u + 2) + 6\sqrt{e_u}\sqrt{e_v}\sqrt{e_w} )v-u)  \\
& -8\sqrt{e_w}(-2u(1-u)w(- 9u^2 + 10u - 2) + 2u^2(1-u)(2-3u) \\
& -\sqrt{e_u}\sqrt{e_v}(3u + 3w - 14uw + 12u^2w - 4u^2) + 2\sqrt{e_u}\sqrt{e_w}u(9u^2 - 14u + 5) \\
& -6\sqrt{e_v}\sqrt{e_w}u(2u^2 - 3u + 1)), \\
C_{d,29} = & -\sqrt{e_u} (48uw(1-u)(1 - w) + 16\sqrt{e_u}\sqrt{e_w}(u(1-w) + 2w(1-u))).   
\end{align*} 
We have that $C_{d,17}(v-u) $ can be bounded by a constant $C(w-u)$ and by convexity, $C_{d,29}\sqrt{e_u}( \sqrt{e_u}-\sqrt{e_v}) \le C(w-u)$. 
%\begin{align*}
%& 16\sqrt{e_u}\sqrt{e_w}(\sqrt{e_u} + 2\sqrt{e_w} - 3\sqrt{e_w}u - 3\sqrt{e_u}w)\blue{(v-u)^3} \\
%&  - 8\sqrt{e_u}\sqrt{e_w}[-2\sqrt{e_u}(3u + 5w - 9uw - 1)+ \sqrt{e_v}(2u + 4w - 6uw - 1)\\
%& + 2\sqrt{e_w}(9u^2 - 10u + 2) + 6\sqrt{e_u}\sqrt{e_v}\sqrt{e_w} ]\blue{(v-u)^2}  \\
%& -8\sqrt{e_w}[-2u(1-u)w(- 9u^2 + 10u - 2) + 2u^2(1-u)(2-3u) \\
%& -\sqrt{e_u}\sqrt{e_v}(3u + 3w - 14uw + 12u^2w - 4u^2) + 2\sqrt{e_u}\sqrt{e_w}u(9u^2 - 14u + 5) \\
%& -6\sqrt{e_v}\sqrt{e_w}u(2u^2 - 3u + 1)]\blue{(v-u)} \\
%& - \blue{e_u(\sqrt{e_u}-\sqrt{e_v})}[48uw(1-u)(1 - w) + 16\sqrt{e_u}\sqrt{e_w}(u(1-w) + 2w(1-u))] 
%\end{align*} 
Therefore, the leading coefficient $K_6(u,v,w)$ is bounded by $C(w-u)$.\\
%Since $\frac{\sum_i \sum_{j \in G_i; j \neq i} (|G_i|-1)^2(|G_j|-1)}{k^2n^4} < \frac{2kn(kn^2)}{k^2n^4} \precsim O(\frac{1}{n})$ and the coefficient is bounded by $C(w-u)$, the entire term is bounded by $C(w-u)^2$. 
\end{enumerate}

Although we have established that the coefficients $K_1(u,v,w), K_2(u,v,w), \hdots K_6(u,v,w)$ can be bounded, in order for the entire expression to be bounded by $C(w-u)^2$ we need the graph configurations in the numerator and denominator to be bounded by $O(1)$ or $O(1/n)$. Recall that the leading term is the denominator is $(nV_G)^2$. Let $\tilde{d_i} = |G_i| - \frac{2|G|}{n}$, then $V_G = \sum_{i=1}^n \tilde{d_i}^2$. 
The graph configurations in the numerator involve: 
\begin{enumerate} 
\item $k^4n^2$
\item $k^2 n \sum_{i=1}^n |G_i|^2$ \label{itm:config2}
\item $\sum_{i=1}^n |G_i|^4$
\item $k \sum_{i=1}^n |G_i|^3$
\item $x_{14}$
\item $\sum_i \sum_{j \in G_i; j \neq i} (|G_i|-1)^2(|G_j|-1)$
\end{enumerate} 

Let $k = O(n^\alpha), 0 \le \alpha <1$. Suppose the largest (centered) degree $\tilde{d}_i \precsim O(n^\beta)$, where $0\le \beta < 1$. 

We first focus on the second configuration \ref{itm:config2} in the numerator, we have: 
\begin{align*}
\sum_{i=1}^n |G_i|^2 & = \sum_{i=1}^n (\tilde{d}_i + \frac{2|G|}{n})^2   \precsim \sum_{i=1}^n (n^\beta + n^\alpha)^2  \precsim n^{2\beta+1} + n^{2\alpha+1}.
\end{align*} 
Since $k^2n \precsim O(n^{2\alpha+1})$, it follows that the entire expression $kn^2 \sum_{i=1}^n |G_i|^2 \precsim n^{2\beta + 2\alpha + 2} + n^{4\alpha+2}$. 

In the denominator, if $\alpha \le \beta$, then $V_G = \sum_{i=1}^n \tilde{d_i}^2 \succsim n^{2\beta}$, and $(nV_G)^2 \succsim n^{4\beta+2}$. Then the ratio of the numerator \ref{itm:config2} and denominator gives us 
$$\frac{n^{2\alpha+2\beta+2} + n^{4\alpha+2}}{n^{4\beta+2}} \precsim O(1).$$ 
If $\alpha > \beta$, then $k^2n \sum_{i=1}^n |G_i|^2 \precsim n^{4\alpha+2}$. % Then we need $(nV_G)^2 \succsim n^{4\alpha+2}$, which in turns implies $V_G \succsim n^{2\alpha} \asymp O(k^2)$. Under this assumption, the remaining configurations (1), (3) - (6), can be bounded accordingly. Explicitly, 
With the assumption that $V_G \succsim k^2 \asymp n^{2\alpha}$, we have $(nV_G)^2 \succsim n^{4\alpha+2}$. Other terms can be done in a similar way. Notice that: 
\begin{enumerate} 
\item[1.] $k^4n^2 \precsim O(n^{4\alpha+2})$. 
\item[3.] $\sum_{i=1}^n |G_i|^4 = \sum_{i=1}^n (\tilde{d_i}+\frac{2|G|}{n})^4
\precsim \sum_{i=1}^n (n^\beta + n^\alpha)^4 \precsim n^{4\beta+1}+n^{4\alpha+1}.$
%When $\alpha \le \beta$, then $n^{4\beta+1}/(nV_G)^2 \precsim  n^{4\beta+1}/(n^{4\beta+2}) \precsim O(1)$. When $\alpha > \beta$, under the assumption that $V_G$ must be at least $O(k^2$), we have $n^{4\alpha+1}/n^{4\alpha+2} \precsim O(1)$. 
\item[4.] $k\sum_{i=1}^n|G_i|^3 \precsim n^\alpha \sum_{i=1}^n (n^\beta + n^\alpha)^3 \precsim n^{3\beta+\alpha+1} + n^{4\alpha+1}$. 
\item[5.] $x_{14} = \sum_i \sum_{j \neq i} (|G_i \setminus \{j \in G_i\} |) (|G_i \setminus \{j \in G_i\}|-1)(|G_j \setminus \{i \in G_j\}|)(|G_j \setminus \{i \in G_j\}|-1) \precsim \sum_{i=1} |G_i|^2 \sum_{j=1} |G_j|^2 \precsim \sum_{i,j}^n (n^\beta+n^\alpha)^4  \precsim n^{4\beta+2} +n^{4\alpha+2}. $
\item[6.] $\sum_{i=1}\sum_{j\in G_i; j\neq i} (|G_i|-1)^2(|G_j|-1) \precsim \sum_{i=1}\sum_{j\in G_i; j\neq i} |G_i|^2|G_j| \precsim n^{3\beta+1+\alpha}$. $\square$

\end{enumerate} 

Therefore, the ratio of the first 5 configurations can be bounded by $O(1)$ and the 6th configuration can be bounded by $O(1/n)$. To see that the 6th configuration can be bounded by $O(1/n)$, consider that if $\alpha \le \beta$, then $(nV_G)^2 \succsim n^{4\beta+2}$ and the ratio of the numerator and denominator is $\frac{1}{n^{(1+\beta-\alpha)}}$. If $\alpha > \beta$, then $(nV_G)^2 \succsim n^{4\alpha+2}$ and the ratio becomes $\frac{1}{n^{(3(\alpha- \beta)+1)}}.$  Recall that expression for $Z_\text{diff}$ can be expressed as the linear combination of the leading coefficients $K_1(u,v,w), \hdots, K_6(u,v,w)$ multiplied by their respective graph configurations.  We have established that $K_1(u,v,w), \hdots, K_5(u,v,w)$ are bounded by $C(w-u)^2$ and $K_6(u,v,w)$ is bounded by $C(w-u)$. Combining these results, and that we are considering the case that $(w-u)>\frac{1}{n}$, it follows that the expression for $Z_\text{diff}$ can be bounded by $C(w-u)^2$. 
\bibliographystyle{imsart-nameyear.bst} % Style BST file (imsart-number.bst or imsart-nameyear.bst)
\bibliography{tightness_ref.bib}
%% or include bibliography directly:
% \begin{thebibliography}{}
% \bibitem{b1}
% \end{thebibliography}

\end{document}